\documentstyle[fullpage, 11pt]{article}

\def\SBIMSMark#1#2#3{
 \font\SBF=cmss10 at 10 true pt
 \font\SBI=cmssi10 at 10 true pt
 \setbox0=\hbox{\SBF Stony Brook IMS Preprint \##1}
 \setbox2=\hbox to \wd0{\hfil \SBI #2}
 \setbox4=\hbox to \wd0{\hfil \SBI #3}
 \setbox6=\hbox to \wd0{\hss
             \vbox{\hsize=\wd0 \parskip=0pt \baselineskip=10 true pt
                   \copy0 \break%
                   \copy2 \break%
                   \copy4 \break}}
 \dimen0=\ht6   \advance\dimen0 by \vsize \advance\dimen0 by 8 true pt
                \advance\dimen0 by -\pagetotal
 \dimen2=\hsize \advance\dimen2 by .25 true in
     \setbox0=\hbox to 3.1 true in{
                \vbox to \ht6{\hsize=3 true in \parskip=0pt  \noindent  
  {\it  Invent. Math.}~{\bf 122} (1995), 1--33
                \vfill}}
  \ht0=0pt \dp0=0pt
 \ht6=0pt \dp6=0pt
 \setbox8=\vbox to \dimen0{\vfill \hbox to \dimen2{\copy0 \hss \copy6}}
 \ht8=0pt \dp8=0pt \wd8=0pt
 \copy8
 \message{*** Stony Brook IMS Preprint #1, #2 ***}
}

\newcommand{\1}{{{\mathchoice {\rm 1\mskip-4mu l} {\rm 1\mskip-4mu l}
{\rm 1\mskip-4.5mu l} {\rm 1\mskip-5mu l}}}}
\newcommand{\R}{{\bf R}}

\newcommand{\Z}{{\bf Z}}

\newcommand{\p}{{\partial}}
\newcommand{\al}{{\alpha}}

\newcommand{\be}{{\beta}}

\newcommand{\Om}{{\Omega}}
\newcommand{\om}{{\omega}}
\newcommand{\eps}{{\varepsilon}}

\newcommand{\ga}{{\gamma}}
\newcommand{\Ga}{{\Gamma}}

\newcommand{\la}{{\lambda}}

\newcommand{\Ll}{{\cal L}}

\newcommand{\Nn}{{\cal N}}

\newcommand{\Pp}{{\cal P}}
\newcommand{\Qq}{{\cal Q}}

\newcommand{\Zz}{{\cal Z}}

\newcommand{\Ham}{{\rm Ham }}

\newcommand{\MS}{{\medskip}}
\newcommand{\BS}{{\bigskip}}
\newcommand{\NI}{{\noindent}}
\newcommand{\proof}[1]{\noindent{\bf Proof#1:\  }}
\newcommand{\jdef}[1]{{\bf #1}}
\newcommand{\QED}{\hfill$\Box$\medskip}
\newcommand{\Cal}{{\rm Cal\,}}
\newcommand{\minset}{{\rm minset}}
\newcommand{\maxset}{{\rm maxset}}
\newcommand{\Tilde}{\widetilde}
\newtheorem{theorem}{Theorem}[section]
\newtheorem{cor}[theorem]{Corollary}
\newtheorem{corollary}[theorem]{Corollary}
\newtheorem{definition}[theorem]{Definition}
\newtheorem{remark}[theorem]{Remark}
\newtheorem{lemma}[theorem]{Lemma}

\newtheorem{prop}[theorem]{Proposition}
\newtheorem{proposition}[theorem]{Proposition}
\newcommand{\at}{{@}}

\title{Hofer's $L^{\infty}$\/-geometry: \\ energy and stability  of Hamiltonian
flows,\\ part I}
\author{Fran\c{c}ois Lalonde\thanks{Partially supported by NSERC grant 
OGP 0092913
and FCAR grant ER-1199.} \\ Universit\'e du Qu\'ebec \`a Montr\'eal
\\ (flalonde \at math.uqam.ca) \and Dusa
McDuff\thanks{Partially supported by NSF
grant DMS 9103033 and NSF Visiting
Professorship for Women GER 9350075.} \\
State University of New York at Stony Brook
\\  (dusa \at math.sunysb.edu)} 
\date{}
\begin{document}
\maketitle

\SBIMSMark{1995/3a}{February 1995}{}
\thispagestyle{empty}

\BS

\NI
{\bf Abstract}
\medskip

Consider the  group $\Ham^c(M)$ of
compactly supported Hamiltonian
symplectomorphisms of the symplectic manifold $(M, \om)$
with the Hofer $L^{\infty}$-norm. A path in $\Ham^c(M)$ will be called a
geodesic if all sufficiently short pieces of it are local minima for the Hofer
length
functional $\Ll$.  In this paper, we give a necessary  condition for a path $\ga$
to be a geodesic.  We also develop a necessary condition for a geodesic to be
stable, that is, a local minimum for $\Ll$.  This condition is related to the
existence of periodic orbits for the linearization of the path, and so extends
Ustilovsky's work on the second variation formula.  Using it, we  construct a
symplectomorphism of $S^2$ which cannot be reached from the identity by a
shortest path. In later papers in this
series, we will use holomorphic methods to prove the sufficiency of the condition 
given here for the  characterisation
of geodesics as well as the sufficiency of the condition for the stability
of geodesics. We will also investigate  conditions under which 
geodesics are absolutely length-minimizing.

\section{Introduction}

 Let $(M, \om)$ be a symplectic manifold without boundary, and let  $\Ham^c(M)$
be the group of
all compactly supported Hamiltonian symplectomorphisms of
$(M,\om)$.  This is an  
infinite dimensional Lie group, whose tangent spaces equal the space of
compactly supported Hamiltonian vector fields on $M$, or,
equivalently,   the space
$$
C_0^\infty(M;\R)/\{{\rm constants}\}
$$
 of compactly supported functions on $M$,
modulo constants.  
In~\cite{HOF}, Hofer considered the  Finsler pseudo-metric
 arising from the norm 
$$
\|H\| = {\rm Totvar}\, H =  \sup_{x\in M} H(x) - \inf_{x\in M} H(x)
$$
on this Lie algebra. He assigned to each $C^\infty$-path 
$\{\phi_t\}_{t \in [a,b]}$
in $\Ham^c(M)$ with $\phi_0 =\1$ the length
$$
\Ll(\phi_t) = \int_a^b {\rm Totvar}\,H_t dt,
$$
where $H_t\in C^\infty(M;\R)$ is its generating Hamiltonian.\footnote
{
Note that this norm is $L^1$ with respect to time $t$ and $L^\infty$
with respect to space.  Eliashberg and Polterovich show in
\cite{EL} that, although one gets an equivalent norm if one varies the norm in 
the $t$-direction, the norm becomes degenerate and essentially trivial if
$L^\infty$ is changed to $L^p$. 
}
Further,
he defined the pseudo-norm $\|\phi\|$ to be the infimum  of $\Ll(\phi_t)$ over
all $C^\infty$ paths $\{\phi_t\}_{t \in [0,1]}$ from $\1$ to $\phi$. 
(This norm is often called the {\bf energy} of $\phi$.)
Setting the distance $\rho(\psi, \phi)$ between two arbitrary points
equal to $\|\phi\circ\psi^{-1}\|$, he obtained  a bi-invariant pseudo-metric
$\rho$ on   $\Ham^c(M)$.

Hofer showed that $\rho$  is indeed a non-degenerate 
metric when $M$ is Euclidean
space $\R^{2n}$ with its standard symplectic structure. 
 In addition, he
showed that the flow $\{\phi_t^H\}_{t\ge 0}$ of an autonomous Hamiltonian $H$
on $\R^{2n}$ is a geodesic with
respect to this norm, in the sense that all sufficiently
short pieces $\{\phi_t^H\}_{s-\eps \le t \le
s+\eps}$  minimize length.  In fact, the path
$\{\phi_t^H\}_{t\in [a,b]}$ minimizes length provided that none of the
symplectomorphisms $\phi_t^H\circ(\phi_a^H)^{-1}, t\in [a,b],$ have non-trivial
fixed points. An appropriate version of this 
result was recently generalised to more
general flows by Siburg in~\cite{SI}. 
Bialy and Polterovich in \cite{BP} improved that result
by a careful analysis of the bifurcations of the action spectrum. 
These proofs use variational
methods which exploit the linear structure of Euclidean space at
infinity.  Thus, other methods are needed in order to extend these
results  to more general manifolds.

In a previous paper~\cite{LALMCD}, we used global embedding techniques and
$J$-holomor\-phic curves to show that $\rho$
is a non-degenerate metric for all $M$.  In this paper and 
its sequels~\cite{LALM1,LALM2},
we will apply these and other techniques to
investigate the properties of geodesics in
$\Ham^c(M)$ for arbitrary $M$, giving in particular a full 
characterization of geodesics and of their
stability, sufficient conditions for geodesics to be absolutely 
length minimizing,
and other related results.   We define geodesics as paths 
which are {\it local}\,\footnote[1]{Throughout this paper, we use the
word \lq\lq local" to mean local in the path space, not local with respect
to  time. A property which holds locally with respect to time will be said
to hold \lq\lq at each moment". }
 minima for $\Ll$  at each moment. In this paper we present those of our
 results which were inspired by a variational approach and are
proved by a variety of {\it ad hoc} techniques.  In particular, we establish
various necessary conditions for a path to be a geodesic by developing several
direct ways in which to reduce the length of a given path.   
We also construct a
symplectomorphism of $S^2$ which cannot be reached from the identity by a
shortest path.  On the other hand,
any result which asserts that a given path is a local or global minimum for $\Ll$
requires  one to measure some associated capacity 
which cannot be reduced.  Our 
results in this direction require new versions of the
non-squeezing theorem which we develop in~\cite{LALM1,LALM2} 
 using  holomorphic methods.   These will allow us to
give conditions under
which a path is length-minimizing, and to establish the sufficiency of the
necessary conditions presented here for a path to be a
geodesic and to be stable.  These results generalize those
 obtained for the case $M = \R^{2n}$ by Bialy-Polterovich
in~\cite{BP} and by Siburg in~\cite{SI}.

\subsection{Geodesics}

Given points $\phi_0,\phi_1\in \Ham^c(M)$, let $\Pp =
\Pp(\phi_0,\phi_1)$ be the space of all
$C^\infty$ paths $\ga = \{\phi_t\}_{t\in [0,1]}$ from $\phi_0$
to $\phi_1$ with the $C^\infty$-topology.  (Thus two paths $\ga$ and $\ga'$ are
close
 if the associated maps $M\times [0,1] \to M$ are $C^\infty$-close.)  For each
$\ga\in \Pp(\phi_0,\phi_1)$ let $\Pp_\ga$ be the
 path-connected component of
$\Pp(\phi_0,\phi_1)$ containing $\ga$.   A path $\ga =
\{\phi_t\}_{t\in [a,b]}$ is said to be \jdef{regular} if its
tangent vector $\dot \phi_t$ is non-zero  for all $t\in [a,b]$.  Further,
$\ga$ is said to be a \jdef{local minimum} of $\Ll$ if it has a neighbourhood
$\Nn(\ga)$ in $\Pp$ 
such that
$$
\Ll(\ga) \le \Ll(\ga'),\;\mbox{ for all }\ga'\in \Nn(\ga).
$$

\begin{definition}\label{def:geod}\rm   Given an interval
$I\subset \R$,  we will say that a
path $\{\phi_t\}_{t\in I}$ is  a \jdef{geodesic}
if it is regular and if every $s\in I$ has a
closed  neighbourhood
$\Nn(s) = [a_s,b_s]$ in $I$ such that the path 
$\{\phi_{\be(t)}\}_{t\in \Nn(s)}$
 is a local minimum
of $\Ll$, where $\be: \Nn(s) \to [0,1]$ is the linear reparametrization
$\be(t) = (t - a_s)/(b_s-a_s)$.
Such a path will be said to be locally length-minimizing at
each moment.  (Thus \lq\lq moments"
have some duration.) A geodesic
$\{\phi_t\}_{t\in [0,1]}$ is said to be   \jdef{stable} if it is a
local minimum for $\Ll$. Note that the notion of stability depends on the
given endpoints of the path, but not the definition of
geodesics.\end{definition}

\begin{remark}\rm (i) We have restricted to regular paths to make it
impossible for a geodesic to stop and then change
direction.  However, this restriction is not essential: see Remark~\ref{rmk:reg}.  
Of course, any regular path may be parametrized by a multiple of its
arc-length without changing its length.   

\NI
(ii)  The above definition has the virtue
that geodesics exist on all manifolds and have a simple characterization: see
Theorem~\ref{thm:geod}.  One might define geodesics in a stronger sense, requiring
that they be absolutely length-minimizing at each moment, instead of locally 
length-minimizing at each moment. Both definitions have their appeal (and 
they agree in ordinary Riemannian geometry).  Our choice was in the end
dictated by the fact that we were unable to establish that 
geodesics in the stronger sense exist on all $M$, though they do 
exist when $M = \R^{2n}$
by the work of Hofer and Bialy--Polterovich (or ours, see \cite{LALM1}).
Another possibility would be to use a variational definition.
Ustilovsky's work~\cite{Ust} shows that this
works very nicely if one restricts attention to paths which satisfy a certain
non-degeneracy condition but, as we shall see below, it
is somewhat cumbersome otherwise.  \end{remark}

  Because Hofer's norm only takes account of the
maximum and minimum values of $H_t$, it is not surprising that the sets
on which $H_t$ assumes these values are important. For each $t\in I$, we
write
\begin{eqnarray*}
{\rm minset} \, H_t = \{ x\in M: H_t(x) = \min H_t\},\\
{\rm maxset} \, H_t = \{ x\in M: H_t(x) = \max H_t\}.
\end{eqnarray*}
A point $q$ which belongs to
$$
\cap_t {\rm minset}\,H_t \quad\mbox{or}\quad \cap_t {\rm maxset}\,H_t
$$
will be called a {\it fixed} extremum of the Hamiltonian $H_t$ over
the interval $I$ and of the corresponding path $\phi_t$.

Sometimes it
is convenient to consider paths $\{\phi_t\}_{t\in [a,b]}$ which do not
start at the identity.  The  Hamiltonian corresponding to such a
path is defined by the requirement that
$$
\frac {d}{dt}\phi_t (x) = X_{H_t}(\phi_t(x))\;\mbox{ for all }t,
$$
where $X_{H_t}$ is the vector field such that
$$
i(X_{H_t}) \om \,=\, \om(X_{H_t}, \cdot) \, =\, dH_t.
$$
Thus it coincides with the Hamiltonian which generates the path
$\{\phi_t\circ\phi_a^{-1}\}$.

Our first theorem characterizes geodesics.

\begin{theorem}\label{thm:geod}  A path
$\{\phi_t\}_{t\in I}$ is a geodesic
if and only if
its generating Hamiltonian has at least one fixed minimum
and one fixed maximum at each moment.  Thus, each $s\in I$ has
a neighbourhood $\Nn_s  \subset I$ such that the Hamiltonian which
generates the path $\phi_t, t\in \Nn_s$, has at least one fixed minimum
and one fixed maximum.
\end{theorem}

We  prove here that
this condition is necessary, postponing to~\cite{LALM1} its sufficiency.  In fact, in
\S2 we describe a simple procedure which shortens every path which does not
have a fixed minimum and maximum. 
The proof that the given
condition is sufficient is more delicate, and relies on a local version of
the non-squeezing theorem  for $J$-holomorphic curves.
This result is already known for the case $M = \R^{2n}$ by the work of 
Bialy--Polterovich~\cite{BP}.\footnote
{
They use rather different terminology, calling
paths with at least one fixed minimum and one fixed maximum  \lq\lq
quasi-autonomous" 
 and paths with a fixed minimum and maximum at
each moment are called \lq\lq locally quasi-autonomous". 
}
It is also proved by Ustilovsky~\cite{Ust} for  paths on an
arbitrary manifold  under the hypothesis
that there is only one fixed minimum $p$ and one fixed maximum $P$ and that 
the Hamiltonian is non-degenerate at these points $p,P$ at all times.

This  characterization of geodesics
 implies that they are not at all unique:  if $\{\phi_t\}$
is a stable geodesic, any path of the form $\{\psi_t\circ\phi_t\}$ will also be
a geodesic
of the same length, provided that the support of $\{\psi_t\}$ 
is disjoint from at least
one pair of fixed extrema $\{p, P\}$, and that $\Ll(\psi_t)$  is
sufficiently small.  Thus we have:

\begin{corollary}\label{profusion}  Given any isotopy 
$\phi_t, 0 \leq t \leq  1,$
there exist an
infinite number of non trivial deformations having the same length. More
precisely, there exists an infinite number of smooth
$1$-parameter deformations $\phi_{t,s}$ such that
\begin{description}
\item[(1)] $\phi_{t,0} = \phi_t$
\item[(2)] $\phi_{0,s} = \1$ and  $\phi_{1,s} = \phi_1$ for all $s$
\item[(3)] for at least one $s$, the isotopy $\phi_{t \in [0,1],s}$ is distinct
from $\phi_{t \in [0,1]}$ and
\item[(4)] for all $s$, $\phi_{t \in [0,1],s}$ has same length as
$\phi_{t \in [0,1]}$.
\end{description}
    In particular, a shortest path or a stable geodesic is never unique.
\end{corollary}

\begin{remark} \rm As Weinstein points out, such non-uniqueness occurs on
a Finsler manifold whenever the unit ball in the tangent space has flat pieces
in its boundary.  A good example to consider is $\R^2$ with the metric whose
unit ball is the unit square $\{(x,y): |x|,|y| \le 1\}$.  
Here, {\it any} smooth
path $(x(t), y(t))$ from $(0,0)$ to $(1,0)$ such that
$$
x'(t) > |y'(t)|
$$
is a geodesic. 
\end{remark}

\subsection{Stability: necessary conditions}

Consider a path $\ga= \{\phi_t \}_{t\in [0,1]}$ with $\phi_0 = \1$.
Suppose that $q$ is a fixed extremum of the Hamiltonian $\{H_t\}_{t\in
[0,1]}$ and   consider the
linearizations
$$
L_t = d\phi_t(q): T_q(M)\to T_q(M)
$$
of the $\phi_t$ at $q$.  Clearly, this is  the symplectic isotopy
generated by the Hessian of $H_t$ at $q$.  What turns out to be crucial
for the stability of $\ga$ is the time at which non-trivial closed
orbits of the $L_t$ appear. If, for every $x\in T_q(M)$ and every
$t'\in (0,T)$, the only trajectories $\al(t) = L_t(x), 0\le t\le t'$, with
$x = L_0(x) = L_{t'}(x)$ are single points, we will say that the linearized
flow at $q$ has no non-trivial  closed trajectories in the time interval 
$(0,T)$.

We first state a necessary condition for stability.

\begin{theorem}\label{thm:stabnec}  Suppose that $\ga$ is a stable
geodesic.  Then it has at least one fixed maximum and one fixed minimum.
Further, if
 $\dim(M) = 2$, there is at least one fixed maximum and  one
fixed minimum at which the linearized flow   has no non-trivial
closed trajectory in the open interval $(0,1)$;  a similar statement
holds for arbitrary $M$ provided that the set of  fixed
extrema of $\ga$ is finite.
\end{theorem}

The first statement follows immediately from the curve-shortening procedure of
Proposition~\ref{prop:curvesh} which reduces the length of every path which does
not have at least one  fixed maximum and minimum.  The second statement is
proved by an explicit construction which shows how to use a closed trajectory
$\al$ of the linearized flow at $q$ to shorten $\ga$.  To do this, one composes
$\ga$ with a  scrubbing motion which moves the points in $M$ lying near $q$
around (the exponential of) the loop $\al$.  Intuitively, in the presence of a closed
trajectory $\al$ at $q$ it costs extra energy to keep $q$ fixed, and one can reduce
the energy needed to get to the endpoint  $\phi_1$ by following $\al$.
The details are in \S~\ref{sec:stab}.

In fact, this was already proved by
Ustilovsky in~\cite{Ust}  under the nondegeneracy assumptions mentioned
before, and our proof uses essentially the same method, but involves more
delicate estimates.   The point is that this nondegeneracy hypothesis on
$\ga$ ensures that the second variation of $\Ll$ at $\ga$ is a well-behaved
functional, and Ustilovsky uses it to prove not only the necessity of the above
condition, but also its sufficiency. 
We  establish the sufficiency of this condition in the general case
in ~\cite{LALM2}.

The above necessary condition places severe restrictions on symplectomorphisms
which are the endpoints of stable geodesics from the identity.  Combining this
with calculations of the  Calabi  invariant of various related symplectomorphisms,
we show:

\begin{prop}\label{prop:nogeod} There is a symplectomorphism $\phi$ of $S^2$ 
which
is not the endpoint of any stable geodesic from the identity.  A fortiori, 
there
is no shortest path from the identity to $\phi$.  \end{prop}

This map $\phi$ is generated by a Hamiltonian of the form
$H(x,y,z)$$ = h(z)$, and so rotates the parallels of the sphere by varying
amounts.

 \subsection{Variational definition
of geodesic}

      Another approach to defining geodesics is to use
a variational definition, looking at paths which are critical points
of the length functional $\Ll$.  In this section we discuss the relationship
between the definition which we have chosen and the variational one.

Observe first that the tangent space, $T_{\ga}\Pp$, to the path
space $\Pp$ at $\ga =
\{\phi_t\}_{t\in[0,1]}$ consists of  smooth
families of functions $G_t,
0 \le t \le 1,$ such that $G_0 = G_1 = 0$.\footnote[1]{Note that when $M$
is non-compact each tangent vector in $T_{\1}Ham^c(M)$ has a unique
representation by a function $G$.  To recover this uniqueness in the
compact case, we  normalise  the function $G$ by requiring that  $\int_M
G\,\om^n = 0$.}
 Further, the tangent  vector $\{G_t\}$
exponentiates to the path $\ga_\eps, |\eps|\le \eps_0,$ in 
$\Ham^c(M)$  given
by   
$$
 \ga_\eps =
\{\phi_{\eps G_{t}} \circ \phi_t\}_{t\in [0,1]},
$$
Here, for each fixed $t$,
$\phi_{\eps G_{t}}$ is the time-$1$ flow of the function $\eps G_{t}$,
and  $\circ$ denotes the usual  composition  of maps.

The following definition takes into account the fact that $\Ll$ is not
differentiable everywhere.  Observe that we do not make a statement about arbitrary
deformations, but only those which arise from exponentiating a vector field as described
above.
 
\begin{definition}\label{def:Lgeod}\rm
A path $\ga = \{\phi_t,\}_{ 0 \leq t \leq 1}$ generated by a
Hamiltonian $H_t$ is said to be \jdef{$\Ll$-critical}
if,  for every tangent vector field $\{G_t\}$, the (not necessarily smooth)
real valued function  $\Ll(\ga_{\eps})$ of the variable $\eps$ is
bounded below on some  neighbourhood of $\eps = 0$ by a  smooth
function whose value at $\eps =0$ is $\Ll(\ga)$ and first derivative
at $\eps=0$ vanishes.  Further, $\ga$ is said to be a smooth point if, for
all tangent vector fields $\{G_t\}$, the
function $\eps \mapsto \Ll(\ga_\eps)$ is differentiable at $\eps = 0$.
\end{definition}

\begin{theorem}\label{thm:Lgeod} A path
$\phi_t, \, 0 \leq t \leq 1,$ is $\Ll$-critical
if and only if its generating Hamiltonian has at least one fixed minimum
and one fixed maximum.
\end{theorem}

Comparing this with Theorem~\ref{thm:geod}, we see that any
$\Ll$-critical path is a geodesic, and that, although a geodesic need not
be an $\Ll$-critical path, it is an $\Ll$-critical path at each
moment. This is in marked contrast with the situation in Riemannian
geometry, where the variational notion of geodesic does not depend on the
interval of time considered.  A path $\ga$
is a Riemannian geodesic exactly when its covariant derivative
vanishes at each time, which implies, of course, that the restriction of
the path to any subinterval, no matter how long,  is also a critical point of the length
functional (on the space of paths with fixed endpoints).
\MS

The next result gives a necessary condition for a $\Ll$-critical path to be 
a smooth point of $\Ll$.

\begin{prop}\label{prop:smooth-point} An isotopy $\phi_t$ generated by
$H_t, 0 \leq t \leq 1,$ with at least
one fixed minimum and maximum,
is a smooth point of the length
functional $\Ll$ only if there exist a fixed minimum $p$ and fixed
maximum $P$ such that
$$
\cap_{t\in [0,1]} {\rm minset} \, H_t \; = \; \{p\} \quad {\rm and}
\quad   \cap_{t\in [0,1]} {\rm maxset} \, H_t \; = \; \{P\}
$$
and such that ${\rm minset} \, H_t = \{p\}$ and ${\rm maxset} \, H_t =
\{P\}$  holds for all $t \in [0,1]$ except on a subset of measure $0$.
\end{prop}

At smooth  points (or more generally at points which satisfy the
hypothesis of continuity, see \S3), Theorem~\ref{thm:Lgeod} follows
directly from the first variation formula of Ustilovsky.  We give the
general proof in \S3.2. Since a direct consequence of Ustilovsky's work
is that, conversely, a path which satisfies the conditions in
Proposition~\ref{prop:smooth-point}
is a smooth point of $\Ll$ provided that each $H_t$ is non-degenerate at both
$p$ and $P$, one sees that the above proposition is close to being sharp.
\MS

\subsection{Organization of the paper}

This paper is organized as follows.  In \S2 we discuss various curve-shortening
techniques and use them to prove the necessary  condition in Theorem~\ref{thm:geod}.
\S3 discusses the first  variation formula for $\Ll$ and proves Theorem~\ref{thm:Lgeod}.
\S4 starts with a discussion of the second variation formula and then proves
the necessary condition for stability in Theorem~\ref{thm:stabnec}. The proof
involves a considerable amount of calculation.   In \S5.1 we apply this theorem to
construct  a symplectomorphism of $S^2$ which cannot be reached by a shortest path
from the identity.    The ideas
in \S2 and \S5 are elementary, and the proofs can be read independently of
everything else in  the paper.
\smallskip

The authors wish to thank Polterovich for some illuminating conversations.
\smallskip

\section{Curve-shortening procedures}

The main aim of this section is to prove the necessity of the condition stated in 
Theorem~\ref{thm:geod} for a path to be a geodesic.  Thus we have to prove
that a path $\{\phi_t\}_{t\in I}$ is
a local minimum for $\Ll$ at each moment  only if its generating
Hamiltonian has at least one fixed maximum and one fixed minimum at each 
moment.  Clearly, this is an immediate consequence of the next proposition.

\begin{prop}\label{prop:curvesh} Suppose that the generating Hamiltonian for
the path $\ga = \{\phi_t\}_{  0 \leq t \leq 1}$ does not have at least one
fixed minimum and one fixed maximum.  Then there is a deformation $\ga_s,
s\ge 0,$ of $\ga = \ga_0$ in $\Pp = \Pp(\phi_0,\phi_1)$ such that
$$
\Ll(\ga_s) < \Ll(\ga),
$$
for all $s > 0$.  In particular, $\ga$ is not a local minimum for $\Ll$.
 \end{prop}

\proof{} By compactness there is a finite set of $t$, say
$t_0 < t_1 < \dots < t_k$ such that
$$
\cap_j {\rm maxset} \, H_{t_j} = \emptyset.
$$
Write $X_j = {\rm maxset} \, H_{t_j}$. Thus,
for some $\nu > 0$
$$
N_{2\nu}(X_0) \; \;\subset \;\; \cup_{j\ge 1} (M - X_j),
$$
where $N_\nu(X)$ denotes the $\nu$-neighbourhood of $X\subset M$ with
respect to some Riemannian metric on $M$. Let $\{\be_j\}$ be a partition of
unity subordinate to the covering
$$
M - N_\nu(X_0), M-X_1,\dots, M-X_k,
$$
and choose
 $\delta > 0$ so that
$$
X_0 \;\; \subset \; \;\cup_{j \ge 1} (\be_j^{-1}([\delta,1]).
$$
For $j \ge 1$, let $K_j$ be a function with support in
$\be_j^{-1}([\delta/2, 1])$ such that
\medskip

\noindent
$\bullet$ $K_j \le 0$,

\noindent
$\bullet$  $K_j$ is constant and $< 0$ on $(\be_j^{-1}([\delta,1])$,

\noindent
$\bullet$  ${\rm supp} \, (K_j) \subset {\rm supp} \, (\be_j)$.
\medskip

Now let $\psi_t^j$ be the time-$t$ flow of $K_j$, and, given $\eps > 0$,
define $\Psi_t^\eps$ as follows:
\begin{description}
\item[(i)]    $\Psi_t^\eps = \1$ for $t < t_0 - \eps $ and then flows
along
$\psi_s^1\circ\dots\circ\psi_s^k$, where $s = t-t_0$, until $t =
t_0+\eps$.

\item[(ii)]  $\Psi_t^\eps$ remains unchanged (its time derivative
is $0$) when $|t - t_j| > \eps$ for all $j$.

\item[(iii)]  When $|t - t_j| \le \eps$, $\Psi_t^\eps$ has the form
$$
(\psi_s^j)^{-1}\Psi_{t_j-\eps}^\eps,\;\mbox{ where }\; s = t - t_j.
$$
Thus as one passes $t_j$ one undoes the $j$th  perturbation.
\end{description}

We claim that $\phi_t^\eps = \Psi_t^\eps\circ
\phi_t$ satisfies the requirements.  Firstly, it is easy to check that
$\Psi_1^\eps = \1$ for all $\eps$.  Further, by (i) and the choice of the
$K_j$, the maximum value of the Hamiltonian  for $
\phi_t^\eps$ is definitely less than that of $\phi_t$ when $|t - t_0 |< \eps$
and $\eps$ is sufficiently small. To see this, note that the Hamiltonian for
the composite $\Psi_t^\eps\circ \phi_t$ is not the sum $H_\Psi + H_\phi$ of
the Hamiltonian for each component but rather is
$$
H_\psi * H_\phi = H_\Psi + H_\phi\circ(\Psi_t^\eps)^{-1}.
$$
However, for small $\eps$, this shifting of the support of $H_\phi$ is
irrelevant in our situation:  $ H_\phi\circ(\Psi_t^\eps)^{-1}$ takes its
maximum on $\Psi_t^\eps({\rm maxset}\, H_t)$ which is contained in $\{x:
H_\Psi(x) < 0\}$ when $|t-t_0| \le \eps$ and  $\eps$ is sufficiently
small.   Thus, there is a constant $c$ which is independent of $\eps$ such
that
$$
\max(H_{\phi_t^\eps}) < \max(H_{\phi_t}) - c
$$
when $|t-t_0|<\eps$ and $\eps$ is sufficiently small.  Similarly,
because the support of $\psi_s^j$ is disjoint from $X_j$ the maximum
 of the Hamiltonian will remain unchanged by the perturbations
described in (iii) for small $\eps$.  It follows that
$$
\Ll(\{\phi_t^\eps\}) < \Ll(\{\phi_t\}) - 2c\eps
$$
as required. \QED

     Note that, as in the above proof, we can compose $H_t$ during the time
interval  $[t - \eps, t + \eps]$ with functions $G_j$ having support disjoint
from the extrema of $H_t$. This proves the non-uniqueness result stated in
 Corollary~\ref{profusion}.

Our next result is a curve-shortening procedure, similar to Sikorav's trick
\cite{S},
which applies to paths with
fixed extrema at which a lot of energy is concentrated.  It gives conditions
under which $\Ll(\ga)$ is not minimal.  Recall that the displacement energy
(or disjunction energy) $e(Z)$
of a subset $Z$ of $M$ is defined by $$
e(Z) = \inf  \{\|\phi\|\,:\, \phi(Z)\cap Z = \emptyset\}.
$$

\begin{proposition}\label{prop:curvsh2} Let $\phi_t$ be a path from $\1$ to
$\phi$  generated by the
Hamiltonian $H_t$ normalised so that 
$\min H_t = 0$ for all $t$, and suppose that there
is $c  > 0$ such that the displacement energy of the set
$$
Z = Z_c = \{ x : H_t(x) \leq c, \;\, \mbox{for some}\;\; t \in [0,1] \}
$$
is less than $c/4$. Assume further that 
$$
\max_{M} H_t \; > \; c/2 + \max_{Z} H_t \quad \mbox{for all} \; \; t.
$$
Then the path $\phi_{t \in [0,1]}$ is not length-minimizing.
\end{proposition}

\noindent
{\bf Proof:} $\; \;$ Let $F:M \to [0,c/2]$ be an autonomous
nonnegative Hamiltonian which equals $c/2$ on $Z_{c/2}$ and has
support in $Z_c$. More precisely, one
may take a set $Z'$ in the interior of $Z_c - Z_{c/2}$, with $\p Z'$
smooth, such that a collar neighbourhood $\p Z' \times [-1,1]$ embeds
in ${\rm Int} (Z_c - Z_{c/2})$.  Let $s$ be the normal coordinate of
the collar chosen so that the points where $s=-1$ are closest to
 $Z_{c/2}$, and define $F(x) = f(s)$ where $f$ decreases from $c/2$ to
$0$ in the interval $[-1/2,1/2]$. 

     By hypothesis, there exists a symplectic diffeomorphism $\tau$ of $M$
of norm less than $c/4$ which disjoins $Z_c$ from itself. Let $\al_t$ be
the Hamiltonian isotopy generated by $-F$, and $\be_t$ the one generated
by $(H_t \circ \al_t) + F$. Set $\al = \al_1$ and $\be = \be_1$. 
Then the path $\al_t \circ \be_t, \, t \in [0,1]$,
where the composition is timewise, is generated by 
$$
  -F + ((H_t \circ \al_t) + F) \circ \al_t^{-1} = H_t.
$$
Thus
\begin{eqnarray*}
\| \phi \| & = & \| \al \circ \be \|  \\
         &=& \| \tau^{-1}(\tau \al \tau^{-1}) \tau \be \| \\
         &=& \| \, \tau^{-1} \left( (\tau \al \tau^{-1}) \be \right) \tau 
               [\tau^{-1}, \be^{-1}] \; \| \\
         &\leq& \| \; [\tau^{-1}, \be^{-1}] \; \| + 
\| (\tau \al \tau^{-1}) \be \| \\
         &<&   c/2 + \| (\tau \al \tau^{-1}) \be \|.
\end{eqnarray*}
   The last inequality holds because 
$$
\| \tau^{-1} \be^{-1} \tau \be \| \; \leq \; \| \tau^{-1} \| + 
                                       \| \be^{-1} \tau \be \|    
                                     \, = \, 2 \| \tau \| < c/2
$$
since the norm is invariant under conjugation. Now the statement of the theorem
follows at once if we show that
$$
\| (\tau \al \tau^{-1}) \be \| \; \leq \; \Ll(\phi_t) - c/2.
$$
But $(\tau \al \tau^{-1}) \be$ is generated by the Hamiltonian
$$
G_t = -F \circ \tau^{-1} + (H_t \circ \al_t + F) \circ \tau \al_t^{-1} \tau^{-1}
= -F \circ \tau^{-1} + F + H_t \circ [\al_t, \tau].
$$
Now over $Z_c$ each function $G_t$ has minimum at least $c/2$: this is obvious
over $Z_{c/2}$, and it holds over $Z_c - Z_{c/2}$ too because each
function $H_t$ is bounded below by $c/2$ there. Because each $H_t$ is bounded
below by $c$ on $M - Z_c$ and since $\tau $ disjoins $Z_c$ from itself and
$F$ has support inside $Z_c$ with values in $[0,c/2]$, 
it is easy to check that the minimum of each
$G_t$ on $M - Z_c$ is also bounded below by $c/2$. Thus 
$$
\min G_t \, \geq \, c/2  \quad \mbox{for all t}.
$$
The same reasons, and the hypothesis that each $H_t$ reaches its maximum
outside $Z_c$ and satisfies $\max_{M} H_t > c/2 + \max_{Z_c} H_t$, imply
easily that $G_t$ has the same maximum value as $H_t$. This concludes the
proof. \QED 

Note that the shorter path from $\1$ to $\phi$ constructed in the above proof
is not $C^{\infty}$-close to the path $\phi_t$. The proof only shows that the
path $\phi_t$ is not length-minimizing, though it might be a local minimum
of the Hofer length $\Ll$ in the path space. We will discuss  the
local minima of $\Ll$  in \S~4.

   Here is an elementary corollary.  Recall  that the 
 the displacement
energy of a ball of radius $r$ in Euclidean space
 is $\pi r^2$.  It follows from~\cite{LALMCD} that this is essentially true for
balls in any
manifold $M$.

\begin{cor}  Suppose that $H$ is an autonomous Hamiltonian which takes its
minimum value at the single point $p$, and suppose that 
$H(x) - H(p) > 4\pi r^2$ for all $x$ outside a symplectically embedded ball $B$ of
radius $r$ and center $p$.    Suppose further that the displacement
energy of $B$  in
$M$ is $\pi r^2$. Then, provided that $\|H\| > 8\pi r^2$, the flow at time $1$
of $H$ is not length-minimizing. \end{cor}

The above  hypothesis will be satisfied if the Hessian of
$H$ at $p$ is large, but we are still quite far from an optimal result.  For
example, in $\R^2$ the function $\pi r^2$ has closed trajectories at time $1$, and
it is easy to see that a
function which equals $\la\pi r^2$ near $0$ will not generate a minimal geodesic
for any   $\la \ge 1$.  But our result only applies when $\la > 4$.

\begin{remark}\rm  (i)
Proposition~\ref{prop:curvsh2} is relevant to the optical
Hamiltonian flows considered by Bialy--Polterovich in~\cite{BP}.  They are
interested in particular Hamiltonians which take their minimum on an
$n$-dimensional section $Z$ of a cotangent bundle $TX$.  When $Z$ is Lagrangian
they show that the corresponding path is always a minimal geodesic.  The above
result makes clear that the Lagrangian condition is essential.  For if $Z$ is not
Lagrangian, it can always be displaced (in fact the displacement energy is $0$ by
Polterovich~\cite{PDIS}), and so if $H$ grows sharply enough near $Z$ the path
will not be a minimal geodesic.

\NI
(ii) This proposition can also be improved in various ways.  For example, it is
clearly unnecessary to assume that the set $Z$ in the statement of the 
proposition contains $\{ H_t^{-1}([0,c])\}$ 
for all $t$ --
if it contains this set  for $t \in [a,b]$ then we should only
\lq\lq turn on" the flow of $F$ for these $t$ as in 
Proposition~\ref{prop:curvesh},
and make corresponding adjustments to the estimates of energy saved.
\end{remark}

\section{$\Ll$-critical paths}

The main aim of this section is to prove Theorem~\ref{thm:Lgeod} which
characterizes $\Ll$-critical paths.  In order to show the logical development
of ideas, we will begin by discussing the first  variation
formula. This has also been derived in a slightly more restricted
context by Ustilovsky~\cite{Ust}.  Since this formula does not apply to all
paths, but only to those which satisfy the Hypothesis of Continuity 
stated below, it is not essential to any of our proofs.  
However, its form is very suggestive.  

\subsection{The first variation formula}

Given $\ga = \phi_{t\in [0,1]}$, let
$G_{t \in [0,1]}$ be a tangent vector field along  $\ga$ vanishing at 
both
ends $t = 0,1$. For any $\eps \in \R$, set
$$
\ga_{\eps}(t) = \phi_{\eps G_{t}} \circ \phi_t
$$
which is a $1$\/-parameter family of paths with the given ends, where
$\phi_{\eps G_{t}}$ is the time $1$ flow of $\eps G_{t}$. We wish to compute
$$
\frac{d^k}{d\eps^k}\Big|_{\eps = 0} \Ll(\ga_{\eps}) =
\int_0^1 \frac{d^k}{d\eps^k}\Big|_{\eps=0} \| \frac{d}{dt} \ga_{\eps} \|_H dt
$$
for $k = 1, 2$.

\begin{prop}\label{prop:taylor} The Taylor expansion of the vector field
$\frac{d}{dt} \ga_{\eps}$
in powers of $\eps$ up to order $2$ is
$$
\frac{d}{dt} \ga_{\eps} = {\rm symplectic \; gradient \; of} \;
{\Huge [}H_t \, + \, \eps \left( G'_t + \{-G_t,H_t\} \right) $$ \vspace{-5mm}
$$ + \frac{\eps^2}{2}
\left( \{-G_t, G'_t\} + \{-G_t,\{-G_t,H_t\}\} \right)  \, + \, o(\eps^2)
{\Huge ]}.
$$
\end{prop}

\NI
Here the notation $o(\eps^k)$ denotes a term $R$ which decreases faster
than $\eps^k$:
$$
\lim_{\eps \to 0}  \frac{R(\eps)}{\eps^k}  = 0
$$
uniformly with respect to other variables.
\MS

\proof{} 
Let $y \in M$ be any point and put $\bar y = \phi^{-1}_{\eps G_{t_0}}(y)$,
where $\phi_{\eps G_{t_0}}$ is the
time$1$-map of the autonomous Hamiltonian  $\eps G_{t_0}$.  We write
$\phi_t^{t_0}$ for 
the flow at time $t$ of the non autonomous Hamiltonian $\{H_{t_0 + t}\}$ (that is:
we look at the flow of the Hamiltonian $H_t$ starting at time $t_0$.)
The vector $\frac{d}{dt}\big|_{t_0} \ga_{\eps}(y)$ is the derivative at $t=0$
of the composition 
$$
[0,\delta] \stackrel{\alpha}{\to} M \times [0,\delta] \stackrel{F}{\to} M
$$
where $\alpha(t) = (\phi_t^{t_0}(\bar y), t)$ and
$F(x,t) = \phi_{\eps G_{t_0 + t}}(x)$.  Now
$\alpha'(0) = (X_{H_{t_0}}(\bar{y}),1)$ where $X$ denotes the symplectic gradient.
Thus
\begin{eqnarray*}
\frac{d}{dt}\big|_{t_0} \ga_{\eps}(y) & = & dF \big|_{(\bar{y},0)}
(X_{H_{t_0}}(\bar{y}),1)\\
& = &  dF \big|_{(\bar{y},0)} (X_{H_{t_0}}(\bar{y})+e_t)
\end{eqnarray*}
where $e_t$ is the unit tangent vector on the real line. Hence
$$
\frac{d}{dt}\big|_{t_0} \ga_{\eps}(y) = d \phi_{\eps G_{t_0}}
\big|_{\bar{y}} (X_{H_{t_0}}(\bar{y})) + \frac{d}{dt} \big|_{t=0}
\phi_{\eps G_{t_0 + t}}(\bar{y}).
$$
Now the first term of the right hand side is equal to
$$
X_{H_{t_0}}(\bar{y}) + \eps [X_{-G_{t_0}}, X_{H_{t_0}}] + o(\eps),
$$
while the second is
$$
 \frac{d}{dt} X_{\eps G_{t+t_0}} + o(\eps)
         = \eps X_{G'_{t_0}} + o(\eps).
$$
Thus finally
$$
\frac{d}{dt} \ga_{\eps} = {\rm symplectic \; gradient \; of} \;
H_t \, + \, \eps \left( G'_t + \{-G_t,H_t\} \right). 
$$ 

A similar but more elaborate calculation shows that the next
term in the Taylor expansion is
$$
\frac{\eps^2}{2}
\left( \{-G_t, G'_t\} + \{-G_t,\{-G_t,H_t\}\} \right).
$$

\QED

 Let
$$
K_{\eps,t} = H_t + \eps ( G'_t + \{-G_t,H_t\} ) + o(\eps)
$$
be the function appearing in the above Taylor expansion.
 To derive the variation formula, we first
make the following assumption.
\medskip

\noindent
{\bf Hypothesis of continuity} \quad The path $\phi_{t\in[0,1]}$,
satisfies the hypothesis of continuity if, given any tangent vector field
$G_{t \in [0,1]}$, it is possible to make
a choice $p_t(\eps) \in M$ of a point at which the minimum value of
$K_{\eps,t}$
is reached in such a way that $p_t(\eps)$ is a smooth path for small values of
$\eps$ and all $0 \leq t \leq 1$. In this case,
$p_t(0) = p_t$ where $p_t$ is a minimal point of $H_t$. We assume that the same
holds for $P_t(\eps)$ with maximum instead of minimum values.
\medskip

One way to decide when this condition is satisfied is to use the following lemma.

\begin{lemma}\label{easy} Let $H_t, 0 \leq t \leq 1$, be any Hamiltonian
which has a non-degenerate minimum at $p$ for all $t$. Then
given any smooth functions $F_t$ of the form $f_t + o(\eps^0)$ defined on
some neighbourhood $N$ of $p$, there is for some $\eps'> 0$  a
smooth map $p(t, \eps): [0,1] \times [-\eps', \eps'] \to N$ such that 
$p(t,\eps)$ is
the unique minimum of $(H_t + \eps F_t) \mid_N$. Further,
$$
\frac{d}{d\eps}\big|_{\eps=0}  \min (H_t + \eps F_t)
=  f_t(p).
$$
A similar statement holds near a non-degenerate maximum $P$.
\end{lemma}

The proof of this lemma is easy, based on ordinary smooth analysis.  It
immediately implies:

\begin{cor}\label{smooth-contin} Suppose that $H_t$ is non-degenerate in
the sense that there exist two points $p,P$ such that for all $t$
\begin{description}
\item[(i)] ${\rm minset} \, H_t = \{p\}$ and ${\rm maxset} \,
H_t = \{P\}$ and
\item[(ii)] $p, P$ are non-degenerate extrema of $H_t$.
\end{description}
Then the path $\ga$ which it generates satisfies the hypothesis of
continuity. \end{cor}
\MS

\begin{theorem}[First variation formula]\label{1var-form} Suppose that 
  $\ga = \phi_{t
\in [0,1]}$ satisfies
the hypothesis of continuity. Then the first variation is:
$$
\delta^1 \Ll (\{G_t\}) = \frac{d}{d\eps}\Big|_{\eps = 0} \Ll(\ga_{\eps}) =
\int_0^1 \left(\sup_{{\rm maxset} \, H_t} G'_t -
\inf_{{\rm minset} \, H_t} G'_t \right) \, dt .
$$
\end{theorem}
 \proof{} Let us compute the total variation of
$$
K_{\eps,t} = H_t + \eps ( G'_t + \{-G_t,H_t\} ) + o(\eps)
$$
 for small $\eps$. 
Under the hypothesis of continuity we find:
\begin{eqnarray*}
{\rm TotVar}(K_{\eps,t}) & = & \left( H_t +
\eps ( G'_t + \{H_t,G_t\} )  \right)
(P_t(\eps))\\
& &\quad\quad - idem (P_t(\eps) \to p_t(\eps)) \, + \, o(\eps).
\end{eqnarray*}
We will write  $\dot P_t(\eps)$ for the derivative of $P_t(\eps)$
with respect to $\eps$.
Then
\begin{eqnarray*}
\frac{d}{d\eps}{\rm TotVar}(K_{\eps,t}) & = & dH_t(\dot P_t(\eps)) +(G'_t +
\{H_t,G_t\})(P_t(\eps))\\ 
& &\quad\quad-  idem (P_t(\eps) \to p_t(\eps))+ o(\eps^0).
\end{eqnarray*}
Therefore, because $dH_t = 0$ at $P_t(0)$ we find that
$$
\frac{d}{d\eps}{\rm TotVar}(K_{\eps,t}) \big|_{\eps = 0} \; = \;
G'_t(P_{t}(0)) - G'_t(p_t(0))
$$
Integrating over $t$ we get the first variation.
Note that $P_t(0)$ is by definition the limit as $\eps \to 0$
of a point where the function
$$
K_{\eps,t} = H_t + \eps \left( G'_t + \{-G_t, H_t\} \right) + o(\eps)
$$
reaches its
maximum  and by the hypothesis of continuity belongs to $\maxset
\, H_t$.   Since $\{G_t, H_t\}$ vanishes over $\maxset\, H_t$,  $P_t(0)$
must belong to the subset 
$$
{\rm maxset} \, (G'_t \mid_{{\rm maxset} \, H_t}) \; \;
\subset \; \; {\rm maxset} \, H_t ,
$$
and similarly for $p_t(0)$.  Therefore
 the first variation formula becomes:
$$
\delta^1 \Ll (\{G_t\}) = \frac{d}{d\eps}\big|_{\eps = 0} \Ll(\ga_{\eps}) =
\int_0^1 \left(\sup_{{\rm maxset} \, H_t} G'_t -
\inf_{{\rm minset} \, H_t} G'_t \right) \, dt.
$$

\QED

\subsection{$\Ll$-critical paths}

Recall that a path $\ga$ is said to be $\Ll$-critical
if,  for every tangent vector field $\{G_t\}$, the (not necessarily differentiable)
real valued function
$\Ll(\ga_{\eps})$ of the variable $\eps$ is bounded below on some
neighbourhood of $\eps = 0$ by a
smooth function whose value at $0$ is $\Ll(\ga)$ and whose first derivative
at $0$ vanishes.
\MS

\NI
{\bf Proof of Theorem~\ref{thm:Lgeod}}

  We must show that a  path
$\phi_t, \, 0 \leq t \leq 1,$ is $\Ll$-critical
if and only if its generating Hamiltonian has at least one fixed minimum
and one fixed maximum.
Suppose that $p,P$ are fixed minimum and maximum of
$\{H_t\}$.
Let $\bar{H}_t, 0 \leq t \leq 1,$ be a $1$\/-parameter family of functions
$C^{\infty}$-close to
$H_t, 0 \leq t \leq 1,$ which is such that for every $t$:
\MS

\noindent
1) $\bar{H}_t(p) = H_t(p)$ and $\bar{H}_t(P) = H_t(P)$

\noindent
2) ${\rm minset} \, \bar{H}_t = \{p\}$ and ${\rm maxset} \, \bar{H}_t = \{P\}$

\noindent
3) $p,P$ are non-degenerate extrema of $\bar{H}_t$ and,
in some symplectic coordinates
near $p$ or $P$, the $2$\/-jet of $\bar{H}_t$ at $p$ is strictly larger
than the $2$\/-jet of $H_t$ at $p$ and conversely at $P$.
\smallskip

As above, for each $\eps $, let $K_{\eps,t}$ be the Hamiltonian which generates
the path
$$
\ga_{\eps}(t) = \phi_{\eps G_{t}} \circ \phi_t,
$$
and set
$$
\bar{K}_{\eps,t} = K_{\eps,t}  - H_t + \bar H_t.
$$
Then, at each fixed minimum $p$, $\min K_{\eps,t} \leq \min \bar{K}_{\eps,t}$.
But $\bar{K}_{\eps,t}$ is now the sum of a function $\bar{H}_t$ which
is non-degenerate at $p$ and of a smooth function $\eps f_t$ (plus 
terms of order
$o(\eps)$), where
$$
f_t = G_t' + \{-G_t, H_t\}
$$
by Proposition~\ref{prop:taylor}. By Lemma~\ref{easy} above,
$$
\frac{d}{d\eps} (\min \bar{K}_{\eps,t}) \Big|_{\eps = 0} =
f_t(p) = G'_t(p) + \{-G_t, H_t\}(p) = G'_t(p).
$$
 A similar result holds at the fixed maxima $P$. Hence we get:
$$
\Ll(\ga_{\eps}) \geq \int_0^1 {\rm Totvar} (\bar{K}_{\eps,t}) \, dt
$$
where the right hand side is a smooth function of $\eps$ whose value
at $0$ is $\|\phi_t\|_H$
and whose first derivative at $\eps = 0$ is therefore  
$$
\int_0^1  (G'_t(P) - G'_t(p)) \, dt  \; = \; 0.
$$
\smallskip

Conversely, if the set of fixed minima or the set of fixed maxima is
empty, one can easily define a tangent vector field $\{G_t\}$ such that
$\Ll(\ga_{\eps})$ is not bounded from below by a smooth function with
vanishing first derivative.  The proof is an obvious adaptation of that of
Proposition~\ref{prop:curvesh}.  Instead of constructing a loop $\Psi_t^\eps$ such
that $\Ll(\Psi_t^\eps\circ\phi_t) < \Ll(\phi_t)$, we now need to find a family
of functions $G_t'$ such that
$$
\int_0^1\;G_t'(x) \,dt = 0\quad\mbox{for all}\; x,
$$
and
$$
\int_0^1 \left(\sup_{{\rm maxset} \, H_t} G'_t -
\inf_{{\rm minset} \, H_t} G'_t \right) \, dt < 0.
$$
If
there are no fixed maxima, for example, there is a finite set of times, say  
$t_ 0
< t_1 < \dots < t_k$ such that
$$
N_{2\nu}(X_0) \; \;\subset \;\; \cup_{j\ge 1} (M - X_j),
$$
where  $X_j = {\rm maxset} \, H_{t_j}$ as before. Then,
 we can choose a small $\delta > 0$ and functions
$$
G_t = \sum_{i = 1}^k G_{k,t} \le 0,\quad\mbox{for }\; |t-t_0| \le \delta,
$$
with support in $N_{2\nu}(X_0)$ so  that
for $i = 1,\dots, k$ and for $|t- t_i| \le \delta$,
$$
{\rm maxset}(H_t - G_{i,t + t_1 - t_0 } ) =  {\rm maxset}(H_t).
$$
It is easy to check that these $G_t$ satisfy the required conditions.
\QED

\begin{cor} \begin{description}
\item[(i)]  A Hamiltonian $\{H_t\}_{t \in [0,1]}$
has at least one fixed minimum and one fixed maximum  if and only if
$$
\int_0^1 \left(\sup_{\maxset H_t} G'_t - \inf_{\minset H_t} G'_t\right) \, dt \;
\geq \; 0 $$
for all admissible tangent vector fields $\{G_t\}_{t \in [0,1]}$.
\item[(ii)] An isotopy $\phi_t$ generated by $H_t, 0 \leq t \leq 1,$ with at least
one fixed minimum and maximum,
is a smooth point of the length
functional $\Ll$ only if there exist two fixed extrema $p,P$ such that
$$
\cap_{t\in [0,1]} {\rm minset} \, H_t \; = \; \{p\} \quad {\rm and} \quad
\cap_{t\in [0,1]} {\rm maxset} \, H_t \; = \; \{P\}
$$
and such that ${\rm minset} \, H_t = \{p\}$ and ${\rm maxset} \, H_t = \{P\}$
holds for all $t \in [0,1]$ except on a subset of measure $0$.
\end{description}
\end{cor}
\proof{} The proof of (i) follows easily from what is said above.
 As for
(ii), if $\phi_t, 0 \leq t \leq 1,$ is
a smooth point of $\Ll$, it has a first derivative which must then be
$$
\delta^1 \Ll (\{G_t\}) =
\int_0^1 \left(\sup_{{\rm maxset} \, H_t} G'_t -
\inf_{{\rm minset} \, H_t} G'_t \right) \, dt.
$$
In particular, this means that the integral expression above is linear in $\{G_t\}$.
If $p \in \cap_{t\in [0,1]} {\rm minset} \, H_t$ and
 $P \in \cap_{t\in [0,1]} {\rm maxset} \, H_t$
is any choice, then
$$
\int_0^1 \left(\sup_{{\rm maxset} \, H_t} G'_t -
\inf_{{\rm minset} \, H_t} G'_t \right) \, dt \; \geq \;
\int_0^1 (G'_t(P) -  G'_t(p)) \, dt \, = \, 0
$$
for all $\{G_t\}$.  If $H_t, 0 \leq t
\leq 1,$ does not satisfy condition (ii), we constructed in the proof of the
Proposition tangent vector fields  $\{G_t\}$ such that the above 
inequality is strict: but then the same integral evaluated on $\{-G_t\}$ cannot be
negative, and so the left hand side cannot be a linear map (it is a singular non-negative
``conic map''). \QED

\section{Geodesics and stability}\label{sec:stab}

We begin this section by discussing the second variation
formula.  Using this as a guide, we then prove Theorem~\ref{thm:stabnec}
which gives a necessary condition for stability. 

\subsection{The second variation formula}

Let $q$ be a fixed extremum of the path $\ga$ at which the Hessian ${\rm d}^2
H_t$ of $H_t$ is non-degenerate for all $t$, and let $\{G_t\}\in T_\ga\Pp$ be a 
tangent vector to $\ga$.    The second variation  of
$\ga$ at  $q$ when evaluated on $\{G_t\}$ depends only on the loop
$g(t)$ traced out by the gradient $\nabla\,G_t(q)$  of $G_t$ at $q$.  (Note that
$g(0) = g(1) = 0$ because  $G_0 = G_1 \equiv 0$.) We
will choose  symplectic coordinates around $q$ and then identify the tangent
space $T_qM$ with $\R^{2n}$ equipped with its standard symplectic form
$\om_0$ and complex structure $J$.  Here $J x_{2i-1} = x_{2i}$, and $Jx_{2i} =
-Jx_{2i-1}$, so that  $$
\om_0(u,v) =  (J u)\cdot v,
$$
where $\cdot$ denotes the usual dot product.  Then
the symplectic area enclosed by a loop $g$ in $\R^{2n}$ is
$$
{\rm area}\, g = \int_{D_g}\om_0 = \frac 12 \int_0^1 (Jg) \cdot g' \, dt,
$$
where $D_g $ is a $2$-disc with boundary along $g$. We
will write  $$
\langle u, v\rangle_t = (( {\rm d}^2 H_t)^{-1} )u\cdot v
$$
for  the metric induced on $T_qM$
by the inverse of the Hessian ${\rm d}^2H_t$ of $H_t$ at $q$.  
The following theorem is proved by Ustilovsky in~\cite{Ust}, and may also be
derived from the Taylor expansion in Proposition~\ref{prop:taylor}.

\begin{theorem}[Second variation formula]\label{2var-form} Suppose that
$H_{t \in [0,1]}$ has at least one fixed minimum and one fixed maximum.
Suppose further that each fixed extremum of $\{H_t\}$ is a 
non-degenerate critical point of all the functions $H_t, 0 \leq t \leq 1$.
Let $G_{t \in [0,1]}$ be a tangent vector field along $\phi_{t
\in [0,1]}$, and set $g(t) = \nabla G_t(p)$.
 Then the contribution of the fixed minimum $p$
of  $H_{t \in [0,1]}$ to the second variational formula, is
$$
\delta^2 \Ll (\{G_t\})(p) = \int_0^1  \langle g', g'\rangle_t\,dt \; + \; 2\, 
{\rm area}\,(g).
$$
Similarly,
the contribution of the fixed maximum $P$ is
$$
\delta^2 \Ll (\{G_t\})(P)   = \int_0^1  
\langle g', g'\rangle_t  \, dt \; - \; 2 \,{\rm area}\,(g) 
$$
where this time $g(t) = \nabla G_t(P)$.
\end{theorem}

We denote by $\Qq_q$ the quadratic
functional 
$$
\Qq_q (g) =  \int_0^1  \left(\langle g',g'\rangle_t
\pm (Jg) \cdot g' \right) \, dt
$$
on the space of smooth loops $g$ based at the origin in
$T_qM = \R^{2n}$ which appears above.   The analysis of this functional
is an isoperimetric problem relating the area of the loop
to its time-dependent energy defined by the varying metric
$\langle\cdot,\cdot\rangle_t$.  It has been carried out
as part of the development of index theory for positive-definite
periodic linear Hamiltonian systems (see Ekeland~\cite{EK}) 
as well as by Ustilovsky in \cite{Ust}.  The results of
the present section show that there is a very close connection between the
periodic linear theory and the question of stability of geodesics in Hofer
geometry. This will become even more apparent in~\cite{LALM2}.

\begin{theorem}[Ustilovsky,\cite{Ust}]\label{Ustilovsky} Let $\ga$ have fixed
non-degenerate extrema $q=p,P$ as above, and suppose that there is no other
fixed extremum.  Then, the quadratic  functional
$\Qq_q$ is positive definite if the linearized isotopy $d\phi_t$ at $q$, generated
by the $2$-jet $\tilde{H}_t, 0 \leq t \leq 1$ of $H_t$ at $q$, has no non-constant
closed trajectory  $ \al$ in time
$\le 1$. Moreover, if this is the case at both $p$ and $P$, $\ga$ is a stable
geodesic, i.e. it is a local minimum for $\Ll$ on the path space $\Pp(\ga).$ 
Conversely, if such $\al$ does exist in
time less than $1$  at
either $p$ or $P$,  then $\Qq_q$ has  non-vanishing index and  the path $\ga$ is
not a  local minimum of $\Ll$. \end{theorem}
\medskip

This theorem.  
can be proved by looking at the $1$\/-parameter family of functionals
$$
\Qq_{t'}(g) =
\int_0^{t'}  \left(\langle g',g'\rangle_t
\pm (Jg) \cdot g' \right) \, dt, \quad
t' \in (0,1]
$$
defined on the space of closed loops
$g:[0,t'] \to T_qM = \R^{2n}$ based at the origin. Note that these functionals
are quadratic (and therefore have generically  only the zero loop as
critical point) and are invariant by translation and
multiplication by $-1$.\footnote[1]{The symmetry group  will
be  larger if, for instance, all the  metrics $\langle \cdot,\cdot\rangle_t, \,
0 \leq t \leq t'$ are conformally equivalent.}

\begin{lemma}\label{le:Ust}  The loop $g$ belongs to the null space of $
\Qq_{t'}$ if and only if $-Jg$ is the translate of  a closed trajectory
$\al$ of $d\phi_t, \, 0 \leq t \leq t'$. \end{lemma}
\proof{}  Let us suppose that $q$ is a minimum so that 
$$
\Qq_{t'} = \int_0^{t'}  \left(\langle g',g'\rangle_t
+ (Jg) \cdot g' \right) \, dt.
$$
Since we may normalise $H_t$ so that its minimum value $H_t(p)$ is $ 0$, its
$2$-jet $\tilde H_t$ may be written in local symplectic coordinates about $p =
0$ as $$
\tilde H_t(x) = \frac 12\sum B_{ij}(t) x_ix_j = \frac 12 x\cdot B_tx,
$$
for some symmetric matrix $B_t= B_{ij}(t)$.  Then the inner product $\langle
u,v\rangle_t =v\cdot (B_t)^{-1}u$, and the linearized flow $L_t =
d\phi_t$ is generated by the vector field $ -JBx$.\footnote[2]{Recall that our
convention is that the symplectic gradient  $X_H$ satisfies
$\om(X_H,\cdot) = dH$.}

Recall that the null space of a quadratic form $\Qq$ on a vector space $V$ is
defined to be
$$
{\rm null \,}\Qq = \{g: \Qq(g,h) = 0 \; \mbox{for all} \; h\in V\}.
$$
Thus $g\in {\rm null\,}\Qq$ if and only if $g$ is a critical point of
$\Qq$.  Now, 
 \begin{eqnarray*} 
\frac{\p}{\p s}\big|_{s=0}
\Qq_{t'}(g+s\xi) & = & \int_0^{t'}\left(2\xi'\cdot B_t^{-1}g' +  J\xi\cdot g'
+ Jg\cdot \xi'\right) dt\\ 
& = & \int_0^{t'} 2\xi'\cdot(B_t^{-1} g' + Jg) dt.
\end{eqnarray*}
Hence, $g$ is in the null space of $\Qq_{t'}$ if and
only if
$$
B_t^{-1} g' + Jg = const,
$$
or equivalently if
$$
g'(t) = B_t (-Jg(t) + c),\quad 0\le t \le t'.
$$
It follows that $-Jg(t) + c$ is  a closed
trajectory  of the linearized flow $L_t, 0\le t\le t'$. \QED

Intuitively, the idea above is that when $\tilde H$, or equivalently $B$, is
small, the term $\langle\cdot,\cdot\rangle_t$ dominates $\Qq_q$.  When $q$ is
a minimum, the closed orbits of the linearized flow $L_t$ enclose negative
area, which increases as $H$ does.  The two terms exactly balance out when
$Jg$  is an orbit of $L_t, 0\le t\le 1$.  
 When $q$ is a maximum the area
enclosed by the closed orbits of $L_t$ is positive, and similar reasoning
applies.
\smallskip

The next step is to show that the
values $t'$  where the null space of $\Qq_{t}$ is non-trivial are
conjugate values.  In other words, for $t < \min t' $, 
$\Qq_{t}$ is positive definite, and the index of $\Qq_{t}$ increases at the
passage of a conjugate value $t'$ by a quantity equal to the (finite) 
 nullity of $\Qq_{{t'}}$.  One can do this by a Lagrange 
multiplier method, or by using the Jacobi sufficient
condition: see~\cite{EK1,Ust}.  This proves the first statement in
Theorem~\ref{Ustilovsky}.

The other statements are proved by investigating
explicit deformations of $\phi_{t\in [0,1]}$ along the loops $g$ in
$T_qM$.  Let $\bar{\al}:[0, t'] \to T_pM = \R^{2n}$
be a closed trajectory of $L_t, \, 0 \leq t \leq {t'}$ and compose it
with some slowing down  function $f:[0,1] \to [0, t']$ which
is the identity on $[0, {t'} - \eps]$ for $\eps>0$ sufficiently
small and sends $[{t'} - \eps,1]$ onto $[{t'} - \eps,{t'}]$. If
$$
\al(t) = \bar{\al}(f(t))
$$
 denotes this composition, define the loop $g$ by requiring that
$$
g(t) = J(\al(t) - \al(0)).
$$
Thus,  $-Jg$ follows a path
which is, up to translation, the same as the path of a closed trajectory of
the linearised isotopy during the time interval $[0, t']$.  Observe that the
choice of $\bar{\al}$ is not unique: it may be replaced by $\rho\bar{\al}$ for
any non-zero scaling factor $\rho$, positive or negative.  

Given such $g$ we define $G_t$ to be a vector
field supported near $p$  with gradient $\nabla G_t(p) = g(t)$.  The
corresponding deformation $\phi_{\eps G_t} \circ \phi_t$ is the composition
of $\phi_t$ with the time-$1$ map of $\eps G_t$. Thus, up to order $1$ in
$\eps$,  
$$ 
\phi_{\eps G_t} \circ \phi_t (p) =
-\eps Jg(t) = \eps (\al(t) - \al(0)).
$$
Ustilovsky showed that it is possible to choose the vector field $G_t$
in such a way that the energy of this deformation is the sum of the energy
$\Ll(\phi_t)$ of the original path with $\Qq(g)$ (up to terms of order
$\eps^3$).  Therefore, if $\Qq(g) <0$, one can decrease the length of
$\phi_t$, while if $\Qq$ is positive
definite one cannot.

The striking  fact here
is that the deformation which optimally reduces the length is given by
composing
the isotopy $\phi_t$ with a motion that moves $p$ in the {\em same} direction
as does the flow of $\phi_t, 0 \leq t \leq 1,$ round $p$.  Thus, if the
linearised motion at $p$ has a closed orbit, the path $\{\phi_t\}$  uses extra
energy to keep the point $p$ fixed rather than letting it move around $p$ 
in the
direction of this orbit.  In the next section we extend the range of validity
of this result, getting rid of most of the non-degeneracy hypotheses on the
path $\ga = \{\phi_t\}$.

\subsection{Stability of geodesics: necessary condition} 
\label{Geodesics-and-stability}

We use the preceding results  as a
guideline to give a rigorous proof of a  necessary condition for the 
stability of geodesics.
For simplicity, we first consider the case when $M$ has dimension $2$.

\begin{theorem}\label{cond-nec-dim2} Let $H_{t \in [0,1]}$ be any
Hamiltonian defined on
a surface $S$, and $\ga = \phi_t, \, 0 \leq t \leq 1$, the corresponding 
path in
$\Ham^c(S)$. If $\ga$ is a stable geodesic,
$H_t$ has at least one fixed minimum $p$ and one fixed
maximum $P$ at which the differential $d\phi_t$ of the isotopy
has no non-trivial closed trajectory in the time interval $(0,1)$.
Indeed, if this condition fails, there
is a canonical deformation of the path $\ga$ which
reduces $\Ll(\ga)$.
\end{theorem}

\proof{}  We have already seen
that a stable
geodesic must have at least a fixed minimum and a fixed maximum. Assuming that
at all fixed  minima of the family $\{H_t\}$ the differentials have a
closed trajectory of period less than $1$, we construct a deformation of the
path $\ga$ which increases the minimum of all $H_t$ without changing
the maxima, and hence reduces the Hofer length of $\ga$. A similar argument
works for maxima.  In the first step of the proof we show how to avoid the
worst degeneracies of $H_t$.    The heart of the proof is Steps 2 and 3
which construct and analyse the scrubbing motion which reduces the length of
the path, and Lemma~\ref{le:remainder} of Step 4 which handles 
the degeneracies of 
$H_t$ at the fixed extremum.  

\medskip

   Let $p$ be a fixed minimum where $L_t = \{d\phi_t(p)\}_{t \in [0,1]}$ has a
non-trivial closed trajectory in
time less than $1$. Observe first that this implies that $p$ is isolated
among the {\em fixed} extrema of $\{H_t\}$, since the manifold is a surface.
Rescale all functions $H_t$ so that their minimum value $H_t(p)$ is $0$.
Note that because $\{H_t\}_{t \in [0,1]}$ defines a geodesic, no function
$H_t$ can be identically zero. Then let $M = \min_t \max_S H_t > 0$ be the 
minimax of the family. 
\medskip

\noindent
{\bf Step 1.}

Working in local coordinates near $p = 0$, let $A = A(\delta)$ be the annulus
$D(4\delta) - D(\delta/2)$ for some small $\delta > 0$, centered at the origin.

\begin{lemma} There exists a deformation of $\ga = \{\phi_t\}$ to a path
(with the same end points and same length $\Ll(\ga)$) which is generated
by a Hamiltonian which is strictly positive on $A$ for all $t$.
\end {lemma}

\proof{}
Suppose to begin with
 that, for at least one value $t_0$, the $2$\/-jet
$\tilde{H}_{t_0}$
of $H_{t_0}$ at $p$ is non-degenerate. We can assume that $t_0 \in (0,1)$
is an
interior value, and that $\delta, \xi$ are small enough so that $H_t$ is 
strictly
positive on $A$ for all $t \in (t_0 - 2\xi, t_0 + 2\xi)$.
Then let $f:N(A) \to [0,1]$ be a $S^1$\/-invariant function defined on a small
neighbourhood of $A$ and
strictly positive and constant on $A$. As in Proposition~\ref{prop:curvesh},
we consider the path $\{\Psi_t\circ\phi_t\}$, where $\Psi_t$ is generated by
the Hamiltonian $F_t = \la(t) f, \, 0 \leq t \leq 1,$ 
where 
\MS

\NI
$\bullet$ $\be: [0,1] \to (-\eps, \eps)$ has vanishing integral; and 

\NI
$\bullet$ $\be$ is equal to
its minimum on $(t_0-\xi, t_0 + \xi)$, and to its maximum on
$[0,1] - (t_0 - 2\xi, t_0 + 2\xi)$. 
\smallskip

Since  $\Psi_0 = \Psi_1 = \1$, the path $\{\Psi_t\circ\phi_t\}$
has the same end points $\1, \phi_1$ as $\{\phi_t\}$. If $\delta$ is chosen
sufficiently small and $\eps$ is smaller than $M/2$,  the
maximum value of the generating Hamiltonian is unchanged and therefore so
is the length $\Ll$. The new path is  generated by a Hamiltonian, that we
still denote $H_t$, which is the same as before everywhere except on $N(A)$
and is always strictly positive on $A$. 

To obtain the same result when all the $2$-jets $\tilde{H}_t$ are
degenerate, it is enough to show that we can
slightly perturb $\{H_t\}$ so that some $H_t$ is strictly positive on $A$.
 But
since a non-constant closed trajectory exists, there must be at least two rank
$1$ functions $\tilde{H}_{t_1},
\tilde{H}_{t_2}$ with distinct kernels: one can then apply the same kind of
argument but using this time a function $f$ which is equal to two bump
functions covering the two connected components of $K_1 \cap N(A)$, where
$K_1$ is the kernel of $\tilde{H}_{t_1}$. This will transform $H_{t_1}$ into a
function strictly positive over $A$ while reducing slightly some positive
values of $H_{t_2}$.
\QED
\medskip

\begin{remark}\label{rmk:reg}\rm In order to make the last step above work, we
used the fact that $H_t$ is not identically $0$ for any $t$.  This is permissible
because the path $\ga$ was assumed to be a geodesic and hence, 
according to Definition~\ref{def:geod}, must be regular.  However, it is
not necessary to assume regularity here: one can use the same trick as
above to make a regular path of the same length as the given one.  To see
this, choose $t_0$ so that $H_{t_0}$ is not identically zero, and let $f:
S \to [0, 1]$ be a smooth function such that $f$ is $0$ on the set of all
fixed minima of the family $\{H_t\}$ and equal to $1$ out of some 
neighbourhood of this set. Then, for all small $\nu$ and all $t$ in some
neighbourhood $(t_0- 2\xi, t_0 +  2\xi)$ of $t_0$, $\max (H_t - \nu f)$ is
reached on the same set as $\max H_t$ and $\nu f < H_t$ everywhere.
It now suffices to compose
$\phi_t$ with $\Psi_t$ generated by $\la f$, where
$\la:[0,1] \to (-\eps, \eps)$
has vanishing integral, reaches its minimum on $(t_0- \xi, t_0 + \xi)$, and
its maximum on $[0,1] - (t_0- 2\xi, t_0 + 2\xi)$.\QED
\end{remark}

 \noindent
{\bf Step 2. Construction of the scrubbing motion} 

By Step 1,
 $m =$$
\inf_t \inf_A H_t$ is $> 0$. 
Since $p$ is a minimum of $H_t$, the linearized isotopy $L_t = d\phi_t$ 
at $p$ of the Hamiltonian
$\tilde{H}_t$ always rotates in the same direction (clockwise, in fact).
Therefore, our hypothesis implies that it
rotates some ray by more than a full turn, and it follows that 
 there exists
 a closed trajectory  $\al:[0,1] \to \R^2 = T_pS$ 
of  $\la \tilde{H}_t$ for some $\la \in (0,1)$.
\footnote[1]{Here we use the parameter
$\la$ as conjugate value parameter instead of 
$t$. In dimension $2$, this will lead to a simpler and more elegant
theory, since there is then a canonical choice of the loop $\al$.
With time as conjugate parameter, one is forced to take a closed
loop $\al:[0, 1] \to T_pS$ obtained by composing the closed
trajectory $\bar{\al}: [0, t'] \to T_pS$ with some more or less
arbitrary slowing down map $f$ as we described in the last section.}
 We construct an optimal deformation of the path
$\phi_t, \, 0 \leq t \leq 1$,
which increases the minimum of each $H_t \big|_{D(4\delta)}$,
by composing $\phi_t$ with a loop $\psi^{\delta,\rho}_t$ which moves the points
near $p$ round a small loop (this is our {\em scrubbing motion}).

For each $t$, and each  sufficiently small $\delta,\rho$,
consider the symplectic diffeomorphism $\psi^{\delta,\rho}_t$ of
$D(3\delta)$ whose restriction to $D(2\delta)$ is the translation by
$\rho \al_0(t)$ where 
$$
\al_0(t) = \al(t) - c, \quad c = \al(0),
$$
and
which is smoothed to the identity on the annulus $D(3\delta) - D(2\delta)$.

We construct the $\psi^{\delta,\rho}_t$ so that they form a closed path, that
is 
$$
\psi^{\delta,\rho}_1 = \psi^{\delta,\rho}_0 = \1.
$$
Thus each
point of $D(2\delta)$ describes a small loop during this Hamiltonian
isotopy. 

Let $F_t$ be the non autonomous Hamiltonian
 which generates the isotopy $\{\psi^{\delta,\rho}_t\}$.    
Since $\psi^{\delta,\rho}_t (x) = x + \rho\al_0(t)$ on $D(\delta)$, the
function $F_t$ must have the form
$$
F_t(x) = \rho  J\al'(t) \cdot x + z(t), \;\mbox{ for }\; x\in D(\delta).
$$
where $z(t) = F_t(0)$.  We normalize $F_t$ by setting $F_t = 0$ on the boundary
of $D(3\delta)$.

\begin{lemma}\label{z}  $\int_0^1 z(t) dt = {\rm area}\, \rho\al_0$.
\end{lemma}

\proof{} Let $\be:[0,1]\to S$ be a path  from a point $\be(0)\in \p
D(3\delta)$ to $\be(1)= 0$.  Then
\begin{eqnarray*}
z(t) & = & F_t(\be(1)) = \int_0^1 \langle dF_t,\dot \be(s)\rangle ds\\
& = & \int_0^1 \om(X_t, \dot\be(s) ) ds,
\end{eqnarray*}
where $X_t = \dot \psi^{\delta,\rho}_t(\be(s))$ and $\om$ is the standard
symplectic form $\om(u,v) = (Ju) \cdot v$ on $\R^2$. Thus $\int z(t) dt$ is the
total flux through the arc $\be$, that is, the total algebraic amount of
surface area which crosses the fixed arc $\al$ during the whole isotopy. This
flux is not only independent of the choice of $\be$ but may also be computed
by taking any family of time dependent arcs $\be_t$, provided that each
$\be_t$ begins at $\be(0)$ and ends at $\be(1)$, and $\be_0 = \be_1$. (Here we
are using the fact that we are working 
locally in $S$ so that the integral of $\om$
over the sphere formed by the images of the paths $\be_t$ is zero.) Take
$\be_t = \la_t \circ \tilde{\la}_t$ where $\tilde{\la}_t$ is the image of the
fixed arc $\be$ by $\psi^{\delta,\rho}_t$ and $\la_t$ is the straight segment
in $D(\delta)$ from $\psi^{\delta,\rho}_t(0) =\rho \al_0(t)$ to $\{0\}$ 
oriented
that way. The total flux is the sum of that through $\tilde{\la}_t$ and that
through $\la_t$.

Since the former follows the flow of the isotopy, the flux crossing it is zero.
To calculate the flux through the arcs $\la_t$ we use the fact that these
paths are entirely contained in the disc $D(\delta)$ on which 
$\phi^{\rho,\delta}_t$ is  translation by  $\rho \al_0(t)$ with  
constant Hamiltonian vector field $X_t =\rho \al_0'(t)$.
The  flux at time $t$ passing through a moving arc $\la_t$ is  the difference
between the infinitesimal flow which passes through $\la_t$ as if $\la_t$ were
fixed, and the infinitesimal area swept out by $\la_t$.  The latter
contribution integrates over $t$ to give the area enclosed by the loop
$\rho\al_0$, while the former is:
\begin{eqnarray*}
\int_0^1\int_0^1 \om(\rho \al_0'(t), \frac {\p}{\p s}\la_t) ds dt & = & 
\int_0^1\int_0^1 \om(\rho \al_0'(t), -\rho\al_0(t)) ds dt\\
& = & -\int_0^1 \rho^2 \al_0'(t)\cdot J\al_0(t) dt = 2\,{\rm area\,}\rho\al_0.
\end{eqnarray*}
\QED

\NI
{\bf Note}
Since the area of $\rho \al_0(t)$ is negative, this average value of
$F_t(q)$, for any $q \in D(\delta)$, is equal to a negative constant.
Of course
the average value at some points in the annulus $D(3\delta) - D(\delta)$
must be positive since the Calabi invariant of the isotopy is $0$. 
\MS

Now consider  the path
$\psi^{\delta,\rho}_t \circ \phi_t$.  It is
generated by
the Hamiltonian $K_t = F_t + H_t \circ (\psi^{\delta,\rho}_t)^{-1}$.  We write
$H_t = \tilde{H}_t + R_t$ on $D(4\delta)$ for all $t$, where $\tilde H_t$ is
the $2$-jet of $H_t$ at $p$ and $\frac{R_t(x)}{\|x\|^2} \to 0$ when $q \to
0$.  Correspondingly, we set 
$$
\tilde{K}_t =
F_t + \tilde{H}_t \circ (\psi^{\delta,\rho}_t)^{-1}.
$$
By construction $K_t = H_t$ outside $D(4\delta)$.  
\MS

\NI
{\bf Step 3.  Calculation of the minimum of $\tilde K_t$ on the disc
$D(\delta)$.} 

For $x\in D(\delta)$,
\begin{eqnarray*}
\tilde{K}_t(x) &=& F_t(x) + \tilde{H}_t \circ (\psi^{\delta,\rho}_t)^{-1}(x)\\
&=&  z(t) + \rho J\al'(t)\cdot x + \tilde H_t(x -\rho\al_0(t))
\end{eqnarray*}
 is a non-homogeneous polynomial of degree $2$.
We now show that its minimum  is reached at
a {\em critical point} lying {\em inside} $D(\delta)$ even when
$\tilde{H}_t$ has rank $1$.
The reason is that we chose  $\al'$ so that $J\al'$ is parallel to the gradient
of $\tilde{H}_t$, which, as we shall see, implies  that the minimum of
$\tilde{K}_t$ may be computed as if the Hessians ${\rm d}^2\tilde {H}_t: \R^2
\to \R^2$ were invertible for all $t$.

\begin{lemma}\label{lambda} There is a continuous path $p(t)$ in $\R^2$ on which $\tilde
K_t$ assumes its minimum over $\R^2$.   By choosing $\rho$ sufficiently
small, we may assume that $p(t)\in D(\delta)$ for all $t$.  Further,
$$
\int_0^1 \min \tilde{K}_t  =
               \int_0^1 (1 - \la) \la \rho^2\tilde{H}_t(\al_0) dt \; > \, 0.
$$
\end{lemma}
\proof{} We prove the lemma
in dimension $2$, but it clearly holds in any dimension.
As in \S4.1, we will write
$$
\tilde H_t(x) = \frac 12 x\cdot B_t x,
$$
for some matrix $B_t$.  Then, the Hessian ${\rm d}^2 H_t$ is the linear
transformation given by the matrix $ B_t$, and 
the closed trajectory $\al$ of the Hamiltonian flow of $\la
\tilde H_t$ satisfies the equation 
$$
\al' = -\la J B_t \al.
$$
Therefore
$$
d\tilde{K}_t (x) = \rho J\al' + d\tilde{H}_t (x -\rho\al) =
B_t (x - \rho\al_0 + \rho\la\al).
$$
This is $0$ when $x \in \rho\al_0 - \rho\la\al + Ker(\tilde{H}_t) $,
and a smooth choice of critical points is given by
$$
p(t) = \rho\al_0 - \rho\la\al = \rho(1-\la)\al - \rho c.
$$
It is clear
that this is small if $\rho$ is small, and that these critical points are
absolute minima of $\tilde{K}_t$ over $\R^2$. 

Observe that
\begin{eqnarray*}
{\rm area}\, \al_0 & = & \frac 12 \int_0^1  J\al_0\cdot \al_0' dt \\
& = &
- \frac 12  \la\int_0^1\al_0\cdot B_t\al_0 < 0.
\end{eqnarray*}
Therefore, by Lemma~\ref{z},
\begin{eqnarray*}
\int_0^1 \min \tilde{K}_t dt  & = & \int_0^1 \tilde K_t(p(t)) \,dt \\
& = & \int_0^1 (z(t) + \rho^2 J\al_0'\cdot (\al_0 - \la\al)+ \frac 12
\rho^2\la^2 \al_0\cdot B_t\al_0)\, dt\\
& = & (-1 + 2(1-\la) + \la)  \rho^2 \la\int_0^1 \al_0\cdot B_t\al_0\,dt\\
& = & (1-\la) \rho^2\int_0^1 \al_0\cdot B_t\al_0\,dt\\
& = & (1-\la)\la \rho^2\int_0^1 \tilde H_t(\al_0) \,dt 
\end{eqnarray*}
is strictly positive because $\al_0$ is a non-constant trajectory of 
$\la \tilde{H}_t$. \QED

\MS

\NI
{\bf Step 4.  The minimum of $K_t$.}

In this step we show how to arrange that $\min\,K_t$ be strictly positive for
all $t$.  To begin we show that
$\int_0^1 \min_{D(4\delta)} K_t$ is strictly positive.

\begin{lemma}\label{le:remainder}  If $\delta$ is sufficiently small, 
we may choose $\rho$ so that
$$
\min_{D(4\delta)} K_t = \min_{D(\delta)} K_t\ge \min_{D(\delta)} \tilde K_t
+\min_{D(2\delta)} R_t 
$$
with
$$
\int_0^1 \left( \min_{D(\delta)} \tilde K_t + 
\min_{D(2\delta)} R_t \right) \, dt 
\; > \; 0.
$$ 
\end{lemma}

\proof{}
 Keeping $\delta$ fixed, and
taking $\rho$ sufficiently small with respect to $m = \min_t\max_A H_t$,
we can insure that the minimum of $K_t\big|_{D(4\delta)}$ is reached inside
$D(\delta)$.  Now $H_t = \tilde H_t + R_t$ on $D(4\delta)$ where
$$
\frac {R_t(x)}{\|x\|^2}\to 0\;\mbox { when }\; x\to 0.
$$
Further $\tilde K_t = F_t + \tilde H_t\circ (\psi^{\delta, \rho}_t)^{-1}$,
where $(\psi^{\delta, \rho}_t)^{-1}(D(\delta)) \subset D(2\delta)$. Thus,
clearly,
$$
\min_{D(\delta)} K_t\ge \min_{D(\delta)} \tilde K_t
+\min_{D(2\delta)} R_t .
$$
We have just seen that $\int \min_{D(\delta)} \tilde K_t$ has the form
$c\rho^2$, where the constant $c$ is independent of $\delta, \rho$.
On the other hand,  $\int \min_{D(2\delta)} R_t = o(\delta^2)$ by the
definition  of $R_t$. Therefore, to prove the second part of the lemma,
it suffices to show that we may choose $\rho = \rho(\delta) $ to be a
linear function of $\delta$.
To check this, consider
 the dependency on $\delta$ of all parameters introduced so far. In Step 1
we introduced a fixed parameter
$\xi$, and parameters $\eps, m$.  These have the form $\eps = {\rm const}
\,  \delta^2$,
$m = {\rm const} \, \delta^2$
since they both only depend on the value of the fixed function $H_t$
(or of the fixed functions
$H_{t_1}, H_{t_2}$) over $A(\delta)$. In Step 2, the functions $F_t$ depend
only on the parameter $\rho = \rho(\delta)$ which determines the size of the
closed orbit.  To insure that the scrubbing motion can be smoothed out to the
identity on $D(3\delta) - D(2\delta)$, one may choose $\rho$ such that
$\max_t \rho\|\al(t)\| \leq \delta/6$, and to be sure that 
$\min_{D(4\delta)} K_t$
is reached on $D(\delta)$, it is enough to choose $\rho$ so that
the minimum over ${D(4\delta)}$ of the linear part of $K_t$ be smaller than
$m/3$, which means that
$4 \delta  \max_t \|\rho \al'\| = 4 \delta \rho \max_t \|\al'\| < m/3
= {\rm const}\, \delta^2$. Thus $\rho(\delta)$ depends linearly on $\delta$, 
as required.
\QED

   We now use  the technique of Proposition~\ref{prop:curvesh} again to deform
the Hamiltonian $K_t, 0 \leq t \leq 1,$ so that $\min_{D(4\delta)} K_t$ is
strictly positive for all $t$.  To do this, 
compose the isotopy with $\psi_t$
generated by the Hamiltonian $F_t, \, 0 \leq t \leq 1,$ defined by
$F_t = \be(t) f$ where
$f:D(4\delta) \to [0,1]$ is a $S^1$\/-invariant bump function equal to
$1$ on $D(3\delta)$ and $0$ near $\partial D(4\delta)$, and where
$\be: [0,1] \to (m_0, m_1)$ has vanishing integral, with
$m_0 = -\max_t \min_{D(4\delta)} K_t$ and $m_1 > -\min_t \min_{D(4\delta)} K_t$.
As before, this composition
has the same end points $\1, \phi_1$, it does not increase the
Hofer length
of the path.  It is now generated by a Hamiltonian, still denoted by
$H_t$, which is the same as before everywhere except on $D(4\delta)$ where each
$H_t$ is now strictly positive.
\MS

\NI
{\bf Step 5.  Completion of the proof of Theorem~\ref{cond-nec-dim2}}.

Repeating the above process near each of the finite number of fixed minima of
$H_t$, we deform $H_t$ to a Hamiltonian $K_t$ with 
$$
\max_x K_t(x) = \max_x H_t(x),\quad \min_{x \in N} K_t(x) > 
\min_{x \in N} H_t(x) = 0,
$$
for all $t$, where $N$ is some neighbourhood of all fixed minima.  
Then, of course, $\{K_t \mid_{S-N} \}_{t\in [0,1]}$ has no fixed minimum,
and Proposition~\ref{prop:curvesh} implies that we can perturb $\{K_t\}_{t \in [0,1]}$
so that their maxima are the same as those of $H_t$, but with minima satisfying
$$
\int_0^1 \min_S K_t > \int_0^1 \min_S H_t.
$$
Thus $\Ll(\{K_t\}) < \Ll(\{H_t\})$.  Further, we may 
clearly choose $\{K_t\}$ to be as close to 
$\{H_t\}$ as we want in the $C^\infty$-topology.
Thus  $\ga$ is not a local minimum of $\Ll$.
\QED

Finally, note that the proof of Theorem~\ref{cond-nec-dim2} shows:

\begin{theorem}\label{cond-nec}  Let $\{H_t\}_{t \in [0,1]}$ be a Hamiltonian
defined on any symplectic manifold $M$, and $\ga = \{\phi_t\}, \, 0 \leq t
\leq 1$ the corresponding isotopy. Assume that each
fixed extremum of $\{H_t\}$ is isolated among the set of
fixed extrema.
If $\ga$ is a stable geodesic,
there exist at least one fixed minimum $p$ and one fixed
maximum $P$  at which the differential of the isotopy
has no non constant closed trajectory in time less than $1$.
\end{theorem}

\proof{}  The proof of Theorem~\ref{cond-nec-dim2}
in the  $2$\/-dimensional case applies directly.
Actually, the hypothesis on dimension has been used only once, namely
to deduce that each fixed extremum is isolated. The only other argument of
the proof which should be treated in a slightly different way is the use
of $t$\/-conjugate values instead of $\la$\/-conjugate values. In arbitrary
dimensions, one cannot derive the existence of a closed trajectory of
$\la \tilde{H}_t, 0 \leq t \leq 1$ from the existence of a closed
trajectory of $\tilde{H}_t, 0 \leq t \leq t'$. Thus, as we indicated
above, the loop $\al$ must be replaced by a closed
loop $[0, 1] \to T_pM$ obtained by composing the closed
trajectory $\bar{\al}: [0, t'] \to T_pM$ with a slowing down map
$f:[0,1] \to [0,  t']$. The rest of the
proof is similar, although the proof of Lemma~\ref{lambda} in Step 3 must be
adapted accordingly.
\QED

\medskip

    This theorem has the following obvious corollary:

\begin{cor}\label{cor:stab} Let $M$ be a compact symplectic manifold, and let 
$\phi \in \Ham(M)$ be
generic in the sense that all its fixed points are isolated.  
 Then, any stable geodesic $\phi_t, \, 0 \leq t \leq 1,$
from the identity to $\phi$
 must have at least two fixed points at which the linearised
isotopy has no non-constant closed trajectory in time less than $1$.
\end{cor}

\section{Symplectomorphisms of $S^2$}

This section is devoted to proving the following result.

\begin{prop}\label{prop:nogeod1} There is
 a symplectomorphism $\phi$ of $S^2$ which
is not the endpoint of any stable geodesic from the identity.  A fortiori, there
is no shortest path from the identity to $\phi$.
\end{prop}

The proof uses properties of the Calabi invariant.  Recall, from~\cite{BAN}
for example, that if $(M,d\la)$ is an exact symplectic manifold, $\Cal$ is a
homomorphism $ \Ham^c(M)\to \R$ defined by:
$$
\Cal(\phi) = \int _{M\times[0,1]} H_t\, \om^n dt,
$$
where $H_t$ is any compactly supported Hamiltonian with time-$1$ map
$\phi$.\footnote[1]{Although this definition does not appear to use the
exactness of $\om$, this is needed to show that $\Cal$ is independent of
the choice of the homotopy class of $\{\phi_t\}$.  For general non-compact
$M$, $\Cal$ is defined on the universal cover of $\Ham^c(M)$.} 
Thus
$$
\Cal(\phi) \le \|\phi\|.
$$

A crucial point  is
that $H_t$ must be compactly supported.  We will see below that if 
$\phi\in \Ham(S^2)$ is
the identity 
near both poles $p_s, p_n$, then the Calabi invariant of $\phi$ considered
as an element of $\Ham^c(S^2 - p_s)$ may be very different from the
corresponding invariant calculated with respect to $\Ham^c(S^2-p_n)$.  It is
this fact which complicates the use of Calabi invariant on $S^2$.\MS

Before starting the construction, we prove the following
easy lemma.

\begin{lemma}~\label{le:area} Let $\{\phi_t\}$ be any isotopy in $\Pp$ with
fixed minimum at $p$ and fixed maximum at $P$, and let $\al$ be a path in
$M$ from $p$ to $P$.  Then $\Ll(\{\phi_t\})$ is the absolute value of the
 area
swept out by $\al$ under the isotopy $\{\phi_t\}$. \end{lemma}
\proof{}  There are several ways to see this.  Here is a geometric argument.
Let $H(x,t)$ be the Hamiltonian which generates $\{\phi_t\}$ and  
consider the surface $S$ in its graph $\Ga_H$ made up of the characteristic
lines starting at the points of $\al$:
$$
S =\{ (\phi_t(\al(u)), H(\phi_t(\al(u)), t), t): u,t\in [0,1]\}.
$$
Then the form $\Om = \om \oplus ds\wedge dt$ vanishes on $S$ since it is a
union of characteristic lines.  Thus
$$
\Ll(\{\phi_t\}) = \int_S ds\wedge dt = -\int_S \om
$$
is (up to sign) the area swept out by $\al$ under the isotopy. \QED

 On the $2$\/-sphere $S$
of radius $1$ centered at the
origin of $\R^3$, take coordinates $\theta: S - \{p_s, p_n\} \to [0,2\pi]$ and
$z: S \to [-1,+1]$, where $\{p_s, p_n\}$ are the south and north  
poles, $\theta (x,y,z)$ is the positive angle of the point $(x,y)$ with respect
to the positive $x$\/-axis, and $z$ is the height coordinate. The symplectic
form is $d\theta \wedge dz$, with total area $A=4\pi$. Thus the Hamiltonian
flow of the function $z$ is the positive rotation
$$
(\theta, z)\mapsto (\theta + t, z).
$$ 

We begin with the
following proposition:

\begin{prop}\label{S2ex} Let $h:\R \to \R$ be a smooth function either strictly
convex everywhere or strictly concave everywhere, with $h'(\pm 1) \notin 2\pi
\Z$, and $\phi$ the time $1$ map of the
Hamiltonian $H = h \circ z$ on $S$.
Then the  length of any stable geodesic
$\psi_t, 0 \leq t \leq 1,$ joining the identity to $\phi$ satisfies
$$
\Ll(\{\psi_t\}) \leq A
$$
\end{prop}

\proof{} \quad Let $\psi_t$ be any stable geodesic from the identity to
$\phi$, generated by a Hamiltonian $K_t$. Let $p,P$ be a fixed minimum
and a fixed maximum of the family $\{K_t\}$, where by Corollary~\ref{cor:stab}
the linearised Hamiltonian isotopy rotates no ray by more than a full turn. Then
$p,P$ belong to $$
{\rm Fix}(\phi) = \{(\theta, z) \mid h'(z) \in 2\pi\Z \; {\rm or} \; z= \pm 1\}
$$
which is the union of a discrete set of parallels.
By Lemma~\ref{le:area}, $\Ll(\{\psi_t\})$ 
is equal to the area swept
out by the curve
$\psi_t(\al(s)), 0 \leq s \leq 1,$ during the time interval
$0 \leq t \leq 1$, where $\al$ is any path from $p$ to
$P$ oriented accordingly. First assume that $p,P$ do not belong to the same
parallel.
Call a path $\al$ from $p$ to $P$ \jdef{admissible}
if it is locally the graph of a function $\theta (z)$: it is a smooth embedded
curve everywhere transversal to the parallels, which can meet the poles $p_n$ or
$p_s$ only at its end-points and only when $\{p_n, p_s\} \cap \{p,P\} \neq
\emptyset$.  Since $h$ is strictly convex or concave,
and $h'(\pm 1) \notin 2\pi \Z$,
the map $d\phi$ at any $q \in {\rm Fix}(\phi)$ has only the tangent
space $T_q{\rm Fix}(\phi)$ as eigenspace. Thus for any
admissible curve $\al$:

a) \, $\al$ intersects $\phi(\al)$ transversally at interior points of $\al$
located on ${\rm Fix}(\phi)$, and all these intersection points have
same sign; and

b) \, $\al(i) = \phi(\al(i)), i=0,1,$ and the tangent vectors are transversal
there.
\MS

   Denote by $\sharp (p,P)$ the algebraic number of interior points of
intersection , which is simply, up to a sign, the number of
parallels in ${\rm Fix}(\phi)$ lying strictly between
$p$ and $P$, thus independent of the choice of the admissible curve.

\begin{lemma} $\sharp (p,P) = 0$.
\end{lemma}

\noindent
{\bf Proof}. \quad  Let $\al$ be an admissible curve from $p$ to $P$. There
is a Hamiltonian conjugation which sends $K_t$ to a Hamiltonian $\hat{K}_t$
on $S$ such that $\hat{p}, \hat{P} = p_s, p_n$, and sends $\al$ to a meridian
$\hat{\al}$. Then $\phi(\al)$ is sent to $\hat{\phi}(\hat{\al})$
which intersects $\hat{\al}$ at $\sharp (p,P)$ interior points of same sign.
Note that the linearised
isotopies $d\hat{\psi}_t$ at $p_s, p_n$ rotate in the positive
$\theta$\/-direction because $p_s, p_n$ are the minimum and maximum
respectively, but no ray turns by more than a full turn.  Further, the tangent
vectors of $\hat{\al}$ and  $\hat{\phi}(\hat{\al})$ at $p_s$ and $p_n$ are
still transversal. Blow-up the sphere at $p_s, p_n$:
 the map $C = ([0,2\pi]/\{0=2\pi\}) \times [-1,1] \to S$ defined by the
coordinates $\theta, z$ admits a unique lifting  of the isotopy
$\hat{\psi}_t(\hat{\al})$ such that $\hat{\psi}_t(\hat{\al}(0))$ is lifted to
$(\theta, -1)$ where $\theta$ is the angle of the tangent vector of
$\hat{\psi}_t(\hat{\al})$ at $\hat{p}=p_s$ (and similarly at $p_n$). Now
lift again to the universal covering $\R \times [-1,1] \to C$ to get
an isotopy $\tau: [0,1] \times [0,1] \to \R \times [-1,1]$ beginning
with Im($\tau(s,t=0)) = \{0\} \times [-1,1]$.
Of course, the condition on the differential
of $\hat{\psi}_t$ at $\hat{p},\hat{P}$ means that $\theta(\tau(i,t)), i=0,1,$
are non-decreasing
functions of $t$ with values in $[0, 2\pi] \subset \R$. But the transversality
of the tangent vectors of $\hat{\al}$ and $\hat{\phi}(\hat{\al})$
at the end points implies that these functions have
values in $(0, 2\pi)$. Since the interior intersection points of
$\tau(s, 1)$ with each of the liftings $\{2\pi k\} \times [-1,1]$ of
$\hat{\al}$ are all transversal and have same sign, the image of
$\tau(s, 1)$ must lie inside $(0, 2\pi) \times [-1,1]$, which means
that $\sharp (p,P) = 0$.
\QED
\medskip

It follows from the proof of this lemma that the area swept out by any 
curve joining $p$ to $P$ is at most
$A$.  By  Lemma~\ref{le:area}, this proves Proposition~\ref{S2ex} when
$p,P$ do not belong to the same parallel.

 If $p,P$ belong to the same parallel, there is no need to
introduce $\sharp (p,P)$: take $\al = \phi(\al)$  a segment of the
parallel to which $p,P$ belong, and the above lifting argument
shows that either $\theta(\tau(0,1)) = \theta(\tau(1,1)) = 0$ or
$\theta(\tau(0,1)) = \theta(\tau(1,1)) = 2\pi$. In the first case, the area
of $\tau$ is $0$, and is $A$ in the second one.  The condition
$h'(\pm 1) \notin 2\pi \Z$ is not necessary, but slightly simplifies the proof:
without it, we would lose the transversality condition of tangent vectors
of $\hat{\al}$ and $\hat{\phi}(\hat{\al})$ at the poles, and we would need
to keep track of the signs to reach the same result.
\QED

\medskip

    Let $h:[-1,1] \to \R$ be any smooth function with $h'(\pm 1) \notin
2\pi (\Z + \frac 12) = \{2 \pi (k + \frac 12) \mid k \in \Z\}$.
Set $\Zz = \{\bar{z} \in (-1,1) \mid h'(\bar{z}) \in
2\pi \Z\}$, and let us denote by $h_{\bar{z}}:[-1,1] \to \R$
the map $h_{\bar{z}}(z) = h(\bar{z}) + \rho (z-\bar{z})$ where $\rho$ equals
$h'(\bar{z})$ if $\bar{z} \neq \pm 1$, and equals $2\pi k$ when
$\bar{z} = \pm 1$ with $k$ the unique integer that minimizes
$|2\pi k - h'(\bar{z})|$. Thus $h_{\bar{z}}$ is the
$1$\/-jet  of $h$ at $\bar{z}$
when $\bar{z} \neq \pm 1$, and is close to the $1$\/-jet when
$\bar{z} = \pm 1$. One should think of $h_{\bar{z}}$ as the correction term
which is needed to make $h$  compactly supported when considered as a function on
$S^2 - \bar{x}$, where $\bar x\in h^{-1}(\bar{z})$.  Of course, $h - h_{\bar
z}$ does not quite have compact support in $S^2 - \bar x$, but its
$1$-jet at $\bar x$ is zero which, as we shall see, means that we can use
it  to calculate the Calabi invariant about $\bar x$ of a slight perturbation
of its time-$1$ map.
Finally, set 
$$
c(h) = - 4\pi +
 \inf_{\bar{z}\in\Zz \cup \{-1,1\}}
\left | \int_{-1}^1 (h - h_{\bar{z}}) \, dz \right |
$$
if the right hand side is positive, and set $c(h) = 0$ otherwise.

The proof of the Proposition now boils down to the following:

\begin{prop} Let $h:[-1,1] \to \R$ be any smooth function with
$h'(\pm 1) \notin
2\pi (\Z + \frac 12)$, and let $\phi$ be the
time $1$ map of $H = h \circ z$. Then any stable geodesic
$\psi_t, 0 \leq t \leq 1,$ joining the identity to $\phi$ satisfies
$$
\Ll(\{\psi_t\}) \geq \frac{c(h)}{2}.
$$
\end{prop}

\noindent
{\bf Proof}. \quad Let $\psi_t$ be a stable geodesic from the identity
to $\phi = \psi$ generated by $K_t, 0 \leq t \leq 1$, with fixed minimum and
fixed maximum $p,P$. We rescale $K_t$ so that $K_t(p) = 0$ for all $t$.
The point $p$ belongs to ${\rm Fix}(\phi) = z^{-1}(\Zz \cup \{-1,1\})$.
We will calculate
in two ways the Calabi invariant about $p$ of a diffeomorphism $\Tilde\phi$
which is very close to $\phi$. 

Suppose first that $z(p) \in \Zz$.
Let
$$
G = (h - h_{\bar{z}}) \circ z
$$
where $\bar{z} = z(p)$. Then $p$ is a critical point of $G$,
$G(p) = 0$ and, because the flows of $h_{\bar z}$ and $h$ commute, $G$ has
time-$1$ map $\phi$. Then let us denote by $\phi_t$ the flow generated by
$G$. Now let  $\be$ be a bump function with
support very near the point $p$, and  $\bar{\phi}_t$ be the isotopy
generated by $\bar{G} = \be G$. (Note that this has very small support.) 
Setting 
$$
 \Tilde G = \left( (1-\be)G \right) \circ \bar{\phi}_t
$$
and denoting by $\Tilde{\phi}_t$ its flow, one easily sees that:
$$
\bar{\phi}_t \circ \Tilde{\phi}_t = \phi_t \quad \mbox{or equivalently} \quad 
\bar G * \Tilde G = G
$$
where $*$ is defined in the proof of Proposition~\ref{prop:curvesh}.  
Thus 
$$
\Tilde{\phi} = \Tilde\phi_1 = \bar{\phi}_1^{-1}\circ \phi
$$
is very close to $\phi$.
Further, because $\Tilde{G} = 0$ near $p$ we may use it to calculate the
Calabi invariant of $\Tilde\phi$
about $p$, that is, the Calabi invariant of $\Tilde\phi$ considered as an
element of $\Ham(S^2 - p)$. We find:
 \begin{eqnarray*}
\Cal_p (\Tilde\phi)& = & \int_t \int_S \Tilde G \om\\
                & = & \int_t \int_S G(\bar{\phi}_t(x)) \om  \quad + \eps_1 \\
       & =  & 2\pi \int_{-1}^1 (h - h_{\bar{z}}) \, dz \quad + \eps_1 \;
\end{eqnarray*}
(In general, in what follows, there will be
various small constants $\eps_i$ which can be made as small as we
want by choosing appropriate bump functions.)

\medskip
   Now let us do the calculation using $K_t$.   We
will add the isotopy $\bar{\phi}_t^{-1}$ to $\psi_t$ (we could tack it on at the
end, that is do $\psi_t$ a fraction faster, and then do $\bar{\phi}_t^{-1}$ quite
quickly) to get an isotopy $\Psi_t$ to $\Tilde \phi$ generated by
$F_t$. Note that $\Psi_t(p) = p$ for all $t$,
and $\Psi_1 = \Tilde{\phi} = \1$ near $p$.

Because $p$ is a minimum of
$K_t$, the rotation of $d\psi_t$ about $p$ is always in the negative
direction.  Also, $\bar{\phi}_t$ is $C^0$-small, and equals the
identity outside the support of $\be$ and on the parallel $z = \bar{z}$
where $G = 0$.  Therefore, it  contributes a total of less than $\pi$ to the
twisting at $p$.  Thus, the isotopy $\Psi_t$ rotates $S$ around $p$
by an angle $\theta_p$ equal either to $0$ or to $-2\pi$.

 Let $\delta$ be a bump
function supported in a little disc centered at $p$, and let $\bar{\Psi}_t$ be the
isotopy generated by $\delta F_t$. As before, let
$\Tilde{\Psi}_t$ be the isotopy generated by
$\left( (1-\delta) F_t \right) \circ \bar{\Psi}_t$. Then 
$\Tilde{\phi} = \Psi_1 = \bar{\Psi}_1 \circ \Tilde{\Psi}_1$. Since all three 
diffeomorphisms fix a neighbourhood of $p$, we can write:
$$
\Cal_p(\Tilde\phi) = \Cal_p(\bar{\Psi}_1) + \Cal_p(\Tilde{\Psi}_1).
$$
Let us begin by computing the first term of the right hand side:
$\bar{\Psi}_1 = \1 $ except on a little
annulus $A$ centered at $p$, whose inner boundary is rotated through angle
$\theta_p$ with respect to its outer boundary.  The isotopy $\bar{\Psi}_t$ fixes
the large disc outside $A$ and moves a small disc near $p$.  However,
to calculate the Calabi invariant of $\bar{\Psi}_1$ about $p$, we must use an
isotopy which fixes a neighbourhood of $p$.  Thus this isotopy must rotate the
large disc outside $A$ through the non-negative angle $-\theta_p$, and therefore,
viewed on the large disc outside $A$, centered at the point antipodal to $p$,
this isotopy rotates the large disc through the non-positive angle $\theta_p$.  It
follows easily that
$$
\Cal_p(\bar{\Psi}_1) = 4 \pi \theta_p + \eps_3.
$$
Therefore
$$
\Cal_p(\Tilde\phi) = \Cal_p(\bar{\Psi}_1) + \Cal_p(\Tilde{\Psi}_1) = 4\pi
\theta_p + c_p +\eps_4
$$
where $c_p = \int_t \int_S (K_t \om)$.
\MS

Thus the two calculations give
$$
\int_{-1}^1 (h-h_{\bar{z}})  \; \leq \; \frac{c_p}{2\pi} \; \leq \;
4\pi + \int_{-1}^1 (h-h_{\bar{z}})
$$
which implies that $\frac{|c_p|}{2\pi} \geq  c(h)$, and therefore that
$$
\Ll(\{\psi_t\}) \; \geq \; \frac{|c_p|}{4\pi} \; \geq \; \frac{c(h)}{2}.
$$

\medskip
   If $p$ is one of the poles, the same argument applies if one
takes $h_{\bar{z}} = h(\bar{z}) + 2\pi k (z - \bar{z})$ where
$\bar{z} = z(p) = \pm 1$ and $k$ is the integer which minimizes
$|2 \pi k - h'(\bar{z})|$. Indeed, the $2$\/-jet of
$G = (h - h(\bar{z}) - 2\pi k (z - \bar{z})) \circ z$ then generates a flow
which rotates the tangent space $T_pS$ by less than $\pi$, and the same
argument goes through. Here again, the hypothesis that
$h'(\pm 1) \notin 2\pi (\Z + \frac 12)$ is not necessary, but slightly
simplifies the definition of $c(h)$.
\QED

\medskip
\NI
{\bf Proof of Proposition~\ref{prop:nogeod1}}

If $h:[-1,1] \to \R$ is  a
strictly convex function with $h'(\pm 1) \notin \pi \Z$, and if the second
derivative $h''(z)$ is large enough (for instance equal to a large constant), then
$\frac{c(h)}{2} > A = 4\pi$, and there cannot exist a stable geodesic joining the
identity to the time $1$ map $\phi$ of $H = h \circ z$.
\QED


\begin{thebibliography}{99999}

\bibitem{BAN}  A.~Banyaga, Sur la structure du groupe des diff\'eomorphismes
qui pr\'eservent une forme symplectique, {\it Comm. Math. Helv.} 
{\bf 53} (1978), 174--227.
 
\bibitem{BP} M.~Bialy and L.~Polterovich, Geodesics
of Hofer's metric on the group of Hamiltonian
diffeomorphisms, preprint, Tel Aviv (1994). 

\bibitem{EK1} I. Ekeland,  An index theory for periodic solutions of
convex Hamiltonian systems, {\it Proc. Symp. Pure Math} {\bf 45} (1986),
395--423.

\bibitem{EK} I. Ekeland,   {\it Convexity Methods in Hamiltonian
Mechanics}, Ergebnisse Math {\bf 19}, Springer-Verlag Berlin (1989).

\bibitem{EL} Y. Eliashberg and L. Polterovich, Biinvariant metrics on the group 
of Hamiltonian diffeomorphisms, {\it Internat. J. Math} {\bf 4} (1993), 727--738.


\bibitem{HOF}  H.~Hofer, Estimates for the energy of a symplectic map,
{\it Comm. Math. Helv.} {\bf 68} (1993), 48--72. 

\bibitem{LALMCD}  F.~Lalonde and D.~McDuff, The Geometry of Symplectic Energy,
       to appear in {\it Annals of Math}.

\bibitem{LALM1}  F.~Lalonde and D.~McDuff,  Hofer's $L^\infty$-geometry: 
energy and stability of Hamiltonian flows II, preprint (1994). 

\bibitem{LALM2}  F.~Lalonde and D.~McDuff, Local non-squeezing theorems and
stability, preprint (1994).

\bibitem{PDIS} L. Polterovich,
      Symplectic displacement energy for Lagrangian submanifolds, 
     {\it Ergodic theory \& dynamical systems} {\bf 13}, (1993), 357--367.

\bibitem{SI} K.~F.~Siburg, New Minimal Geodesics in
the Group of Symplectic Diffeomorphisms, to appear in {\it Calculus of Variations}.

\bibitem{S}  J.~C.~Sikorav,
Syst\`emes hamiltoniens et topologie symplectique,
              ETS Editrice, Pisa, 1990. 

\bibitem{Ust} I. Ustilovsky, Conjugate points on geodesics
of Hofer's metric,
preprint, Tel Aviv, (1994).


\end{thebibliography}
\end{document}